\documentclass[12pt]{article}

\usepackage{amsfonts}

\def\a{\mathfrak{a}}

\def\si{\sigma}
\def\ep{\varepsilon}

\def\vid{\emptyset}

\def\A{\mathcal{A}}

\def\LL{\Lambda}
\def\l{\lambda}
\def\aa{\alpha}

\def\cc{\check}

\def\bb{\backslash}

\def\F{\bold{F}}

\def\mm{\mathcal}

\def\N{\mathbb{N}}

\def\Z{\mathbb{Z}}
\def\R{\mathbb{R}}
\def\C{\mathbb{C}}

\def\fd {\hspace{0.35cm} \raise
-0.5mm\hbox{$\blacksquare$}\\}

\def\qed{{\null\hfill\ \raise3pt\hbox{\framebox[0.1in]{}}}\break\null}
\newtheorem{theo}{Theorem}
\newtheorem{prop}{Proposition}
\newtheorem{lem}{Lemma}
\newtheorem{cor}{Corollary}
\newtheorem{rem}{Remark}
\newtheorem{defi}{Definition}

\def\ste{\par\smallskip\noindent}

\def\dem{ {\em Proof~: \ste }}

\def\beq{\begin{equation}}
\def\eeq{\end{equation}}

\def\bb{\backslash}

\def\A {\cal A}

\newenvironment{res}
               {\begin{equation}\begin{minipage}{0.85\textwidth}}
               {\end{minipage}\end{equation}}
\def\ber{\begin{res}}
\def\eer{\end{res}}

\def\un{ \underline}
\def\mm{\mathcal}

\catcode`\Ž=\active\def Ž{\'e}
\catcode`\=\active\def {\`e}
\catcode`\ˆ=\active\def ˆ{\`a}
\catcode`\=\active\def {\`u}
\catcode`\=\active\def {\^{e}}
\catcode`\"=\active\def "{\^{i}}
\catcode`\™=\active\def ™{\^{o}}
\catcode`\ž=\active\def ž{\^{u}}
\catcode`\=\active\def {\c{c}}
 \def\opto#1{{\smash{\mathop{\longrightarrow}\limits^#1}}}

\title{Neighborhoods at infinity  and the Plancherel formula for a reductive $p$-adic symmetric space}

\author{ Patrick Delorme}

\date{ abreviated title: Plancherel formula for a  reductive$p$-adic symmetric space}

\begin{document}

\maketitle

 Aix Marseille Universit\'e\\
CNRS-IML, FRE 3529\\
163 Avenue de Luminy\\ 13288 Marseille Cedex 09\\  France.\\          
 delorme@iml.univ-mrs.fr
 \\\\
 The author has been supported by the program ANR-BLAN08-2-326851 and by the Institut Universitaire de France  during the elaboration of this work.
 \\\\ {\bf Abstract} Yiannis  Sakellaridis and Akshay Venkathesh have  determined, when the group $G$ is split  and the field $\F$ is of characteristic  zero, the Plancherel formula for  any spherical space $X$ for $G$ modulo the knowledge of the discrete spectrum.  
 
The starting point is the determination of good neighborhoods at infinity of $X/J$, where $J$ is a small  compact open subgroup of $G$.  These neighborhoods are related to ''boundary  degenerations'' of $X$. The proof of their existence is made by using wonderful compactifications.

In this article we  will show the existence of such neighborhoods  assuming that  $\F$ is  of characteristic different from 2 and $X$ is symmetric.
In particular, one does not assume that $G$ is split. Our  main tools are  the Cartan decomposition of Benoist and Oh,   our previous  definition of the constant term  and asymptotic properties   of Eisenstein integrals due to Nathalie Lagier . 

Once the existence of these neighborhoods at infinity of $X$ is  established,  the analog of the work of Sakellaridis and Venkatesh  is straightforward and leads to the Plancherel formula for $X$. 
\\\\{\bf Classification} 22E

\newpage
\section{Introduction}
Let $G$ be the  group of $\F$-points of a reductive group $\un{G}$ defined over the non archimedean local field $\F$.

In a tremendeous work (cf. [SV]), Yiannis  Sakellaridis and Akshay Venkathesh have  determined, when the group $G$ is split  and the field $\F$ is of characteristic  zero, the Plancherel formula for  any spherical space $X$ for $G$ modulo the knowledge of the discrete spectrum.

The starting point is the determination of good neighborhoods at infinity of $X/J$, where $J$ is a small  compact open subgroup of $G$. Notice that $G$ acts on $X$ on the right.  These neighborhoods are related to ''boundary  degenerations'' of $X$. The proof of their existence is made by using wonderful compactifications 

In this article we  will show the existence of such neighborhoods  assuming that  $\F$ is  of characteristic different from 2 and $X$ is symmetric.
In particular, one does not assume that $G$ is split. The main tool is the Cartan decomposition (cf. [BenO]), the definition of the constant term (cf [D])  and asymptotic properties   of Eisenstein integrals due to Nathalie Lagier (cf. [L]). 
The use of Eisenstein integrals to prove results geometric in nature on symmetric spaces goes back to her  work (cf. [L]  Theorem 7).
Notice that our  neighborhoods at infinity are quite explicit in terms of the Cartan decomposition. 

Once the existence of these neighborhoods at infinity of $X$ is  established,  the analog of part 3 of [SV] is straightforward and leads to the Plancherel formula for $X$. Notice that our definition of normalized integrals differs slightly from the one in [SV] section 15.

Let $\si$ be an  involution of $\un{G} $ defined over $F$. Let $H$ be the fixed point group of $\si$ in $G$ and let $X= H\bb G$.
We denote by $X(G)$  the group of unramified characters of $G$ and $X(G)_\si$ be the connected component of $1$ in $\{\chi\in X(G)\vert  \chi\circ \si= \chi^{-1}\}$.

Let  $P$ be a $\si$-parabolic subgroup of $G$ i.e. such that  $P$ and $\si(P)$ are opposed.  Let $M:=P\cap \si(P)$ be the $\si$-stable Levi subgroup of $G$.  Let $U$ ( resp., $U^-$ ) be the  unipotent radical of $P$ (resp., $P^-:= \si(P)$) and let $\delta_P$ be the modulus function of $P$.
We define:
$$H_P= U^- (M\cap H),  X_P= H_P \bb G$$ 
 The space $X_P$ is called a ''boundary degeneration '' of $X=H\bb G$. It is an important object  whose role has been emphasized by  Sakellaridis and Venkatesh.
 
 Let $P_\vid=M_\vid U_\vid $ be a minimal $\si$-parabolic subgroup of $G$. We will assume in this introduction that $P_\vid H$ is the only $(P_\vid ,H)$-open  double coset in $G$. A   split torus is said $\si$-split if  all its elements are antiinvariant by $\si$.
Let $A_\vid$ be the  maximal $\si$-split torus of  the center of $M_\vid$ . Let $A_\vid ^+$ be  the closed  positive chamber in $A_\vid$ for $P_\vid$.  
 The Cartan decomposition asserts (cf. [BenO]):
 $$ G= H A_\vid^+ \Omega, $$
 for some compact subset, $\Omega$,  of $G$. 
 Let $P$ be a $\si$-parabolic sugroup of $G$.
 If $C>0$, let $$A^+_\vid (P, C):=\{ a \in A_\vid^+ \vert \vert \aa (a) \vert _F > C, \aa \>root \> of  \> A_\vid \> in \> U \}.$$
 
 We denote by $\dot{1}$ (resp., $\dot{1}_P $) the image of the neutral element 1 of $G$ in $X$ (resp., $X_P$). 
 The following Theorem  (cf. Theorem \ref{theoct}) is an easy consequence of [D], Proposition 3.14.
\\{\bf Theorem (Constant term map)} \\{\em There is  a unique $G$-equivariant map $c_P:C^{\infty}(X) \to C^{\infty } (X_P)$ with the following property.
For every compact open subgroup $J$ of $G$, there exists $C>0$ such that for all $f \in C^{\infty}(X)$ which is $J$-invariant:
$$(c_Pf)(\dot{1}_P a \omega)= f(\dot{1}  a \omega), a \in A_\vid^+(P,C) , \omega \in \Omega.$$ }
The following theorem (cf. Theorem \ref{theoexp} for its detailed version) was suggested by the work [SV] of Sakellaridis and Venkatesh, who constructed similar maps, in their context,  using wonderful compactifications.
 \\ {\bf Theorem   ($exp_{P, J}$-maps)    } \\ {\em Let $P=MU$ be a  standard $\si$-parabolic subgroup of $G$ i.e. such that $P_\vid \subset P$. Let $J$ be a compact open subgroup of $G$.\\
(i)  There exists $C>0$ such that the correspondence  $\dot{1}xJ \mapsto \dot{1}_P x J$, for $x \in A_\vid^+(P,C) \Omega$, is a well defined  bijective map denoted $exp_{P,J}$
from the  subset  $N_{X,J}(P,C):= \dot{1} A_\vid^+(P,C)\Omega J$ of $X /J$, to the  subset $O':=\dot{1}_P A_\vid^+(P,C)\Omega J$  of $X_P /J$. \\
 (ii) For $J$ small enough,  the map $exp_{P,J}$ preserves volumes. 
\\ (iii) For $f $ any  right $J$-invariant element  of $ C^\infty (H\bb G)$, one has :
$$(c_P) (exp_{P,J} (x))=  f(x), x \in N_{X,J}(P,C)  . $$ }

As said above we need some results of N. Lagier on Eisenstein integrals that we will recall.  Let $P=MU$ be a $\si$-parabolic subgroup of $G$. Let $(\delta,E)$ be   a unitary irreducible  representation  of $M$.  Let $\chi \in X(M)_\si$ and let $\delta_\chi= \delta \otimes \chi$. We denote by  $i^G_P\delta_\chi $  or $\pi_\chi$ the   normalized induced representation and let $V_\chi $ denote its space.  
   
Let $\eta\in E'^{M\cap H}$. Let $\chi\in X(M)_\si$, sufficiently $P$-dominant. There is a canonical $H$-fixed linear form $\xi(P, \delta_\chi, \eta)$ on $V_\chi$,  (cf. [BD]). One defines the Eisenstein integrals on $X$, $E(P, \delta_\chi, \eta, v)\in C^{\infty} (X), v \in V_\chi.$ by:
$$  E(P, \delta_\chi, \eta, v) (\dot{1} g)= \langle \xi(P, \delta_\chi, \eta), \pi_\chi (g) v\rangle, g \in G $$
Let  $A_M$  be the maximal $\si$-split torus of the center of $M$ and let  $\mu_\delta$ be the character of $A_M$ by which $A_M$ acts on $\delta$. 
The following theorem  is due to Nathalie Lagier.  This is the analog of a Langlands lemma on asymptotics of smooth coefficients.
\ber 
{\em One says that the sequence   $(a_n) $ satisfies $(a_n) \to_P\infty$ if   $a_n\in A_M$ and for every root  $\aa$ of $A_M$ in the Lie algebra of $U$, $(\vert \aa (a_n)\vert_F)$ tends to infinity. 
\\ Let us assume $Re (\chi)\delta_P^{-1/2} $ is  $P$-dominant. Then if $(a_n) \to_P \infty$ the following limit exists $$lim_{n\to \infty}  (\chi \delta_P^{-1/2} )(a_n^{-1})\mu_\delta(a_n^{-1}) E(P, \delta_\chi, \eta, v) (\dot{1} a_n)$$
and is equal to 
$$< \eta, (A(P^-, P,\delta_\chi) v )(1)>, $$
where $A(P^-, P,\delta_\chi)$ is  the  (converging) intertwining integral operator.  }\eer 
The theorem admits a  variation  when $(a_n) \to_{Q}\infty$, with $P\subset Q $.
This implies easily: \\ 
{\bf First Key Lemma}\\
{\em Let us  assume  that $(a_n) \to_P \infty $ and  that $(g_n)$ is a sequence in $G$ converging to $g$.  If $\dot{1} a_n g_n = \dot{1} a_n $ for all $n$, then $g \in U^-(M\cap H)$. } 
\\ {\bf Second Key Lemma}
 \\ {\em Let $J$  be a compact open subgroup of $G$. Let $(a_n) \to_P \infty,  (a'_n) \to_{P'} \infty$, for  $P, P'$ $\si$-parabolic subgroups of $G$ and let   $g, g'\in G$. Let us assume $\dot{1} a_n g J= \dot{1}  a'_n g' J, n \in \N$.\\
  Then $P=P'$ and a subsequence of  $({a_n }^{-1} a'_n )$ is bounded. }
\\ {\bf Definition  of $exp_{P,J}$}
 Although we gave a formula for $exp_{P,J}$ it is unclear that it is well defined. Let us sketch the proof that it is well defined. 
If it was not well defined  for all $C>0$, there would  exist  two  standard  $\si$-parabolic subgroups $Q$, $Q' $ of $G$  contained in $P$, and  two sequences  $(a_n) \to_Q \infty,( a'_n)\to_{Q'} \infty $, $u, u'\in G$ such that 
$$\dot{1} a_n uJ =\dot{1} a'_nu'J$$
and $$\dot{1}_P a_n u' J\not= \dot{1}_P a_n u'J.$$
By the  Second Key Lemma,    one sees that $Q=Q'$ and, possibly   extracting a subsequence,  one has  from the First Key Lemma:     $$\dot{1}_Q  a_nu J= \dot{1}_Q   a'_n u'J.$$ A trick  (see  below for an other occurence of this trick)  with the constant term of the characteristic function  of a $J$-coset in $X$  allows to show $\dot{1}_P a_nu J= \dot{1}_P a'_n J$ for $n$ large. A contradiction which proves our claim. \\
{\bf Injectivity of $exp_{P,J}$}
\\One wants to prove that, for $C$ large, if  $x, x' \in N_{X,J}(P,C)$  and $ exp_{P,J}(x)= exp_{P,J}(x')$, then $x=x'$. One introduces the characteristic function $f$ of $ x \subset X$  and  one will use its constant term $c_Pf$. These functions are $J$-invariant and their values on a $J$-coset makes sense.  From the properties of the constant term, if $C $  is large enough one has:
$$(c_Pf)( exp_{P,J}(x))= f(x)=1 .$$  
 But,   from our hypothesis one deduces : 
$$(c_Pf)(exp_{P,J}(x) )= (c_Pf)( exp_{P,J}(x') ).$$ 
Moreover  by the properties of the constant term and because 
$C$ is large, one has: 
$$(c_Pf)( exp_{P,J}(x'))= f(x')=1 .$$ 
This implies that  $f(x')=1$, hence $x=x'$, as wanted. \\

 A  compact open subgroup $J$  of $G$ is said to have a strong $\si$-factorization   for $P_\vid$ if for all $\si$-parabolic subgroup $P=MU$ which contains $P_\vid$ one has:\\
1)$J= J_{U^-} J_M J_U$  for all $\si$-parabolic subgroups, where $J_M= J\cap M,...$, 
\\2) For all $a\in   A_\vid ^+ $, $a^{-1}J_U a  \subset J_{U}$,  $aJ_{U^-} a^{-1} \subset J_{U^-}$.
\\3) $J= J _H J_P$, where $J_H= J \cap H, J_P= J \cap P$.
\\4) $J_M$ satisfies the same properties for $P_\vid \cap M$. 
\\There are arbitrary small compact open subgroups with a strong $\si$-factorization  for $P_\vid$ (cf. Kato-Takano [KT1] if the residual characteristic  is different from 2, [CD] in general and Lemma \ref{strongsi} of this article  for the ''strong'' version).  

A choice of  a $G$-invariant measure on $X$ determines a $G$-invariant measure on $X_P$.  
\\  {\bf Third Key   Lemma}
{\em Let  $J$ be a compact open subgroup with a strong $\si$-factorization  for $P_\vid$. Let   $a \in A_\vid ^+$. Then:
 $$\dot{1}  a J=  \dot{1} a  J_{M} J_{U}, \dot{1}_P  a J=  \dot{1}_P  a  J_{M} J_{U}, $$
 $$ vol_X  (\dot{1} aJ)= vol_{X_P} (\dot{1} _P a J).$$}
 The proof is easy. Moreover one can show  that the identity of volumes is also true for any small enough compact open subgroup of $G$. This implies easily the last  property of $exp_{P, J}$. 
 
 Then, one introduce the restriction $e_P$  of the transpose map of the constant term map  to $C_c^\infty(X_P)$. Following an idea given to us by Joseph Bernstein,  and using a result of Aizebbud, Avni, Gourevitch[AAG] , one shows that its image is contained in $C_c^\infty(X)$ (cf  Theorem \ref{theoe}). This achieves to prove the analog of Theorems 5.1.1 and 5.1.2 of [SV].
 
 Then, as was said before, this allows to prove the Plancherel formula, modulo the discrete spectrum of the $X_P$, by using the same method than [SV], Part 3.
 
 It is natural to ask about spherical varieties for general reductive groups. We think that the key point is the existence or not of a Cartan decomposition. 
 \\{\bf Acknowledgments} I thank warmly Joseph Bernstein for very useful discussions. I thank also Yiannis Sakellaridis who answered to my numerous questions about his work with Akshay  Venkatesh.  I thank also Pascale Harinck and  Omer  Offen for useful  discussions about the case $GL(n, F) \bb GL(n,E)$.
\section{Notations}

\setcounter{equation}{0}
If $E$ is a vector space, $E'$ will denote its dual. If $T:E\to F$ is a linear map between two vector spaces, $T^t$ will denote  its transpose.   If $E$ is real, 
$E_\C$ will denote its complexification. If $G$ is a group,
$g\in G$ and  $X$ is a subset of $G$, $g.X$ will denote $gXg^{-1}$. If $J$ is a subgroup of $G$, $g\in G$
and $(\pi, V)$ is a representation of $J$, $V^J$ will denote the space of
invariant elements of $V$ under $J$ and $(g\pi, gV)$ will denote the representation of $g.J$ on $gV:=V$ defined by: $$(g\pi)(g.x): = \pi(x), x\in J.$$ We will denote by $(\pi', V')$ the contragredient representation of a representation $(\pi,V)$ of  $G$ in the algebraic dual vector space $V'$ of $V$. \\If $V$ is a vector space of vector valued functions  on $G$ which is invariant by right (resp., left) translations, we will denote by $\rho$ (resp., $\l$) the right (resp., left) regular representation of $G$ in $V$.\\ If $G$ is locally compact, $d_lg$ or $dg$ will denote a left invariant Haar measure on $G$ and $\delta_G$ will denote the modulus function.
 
Let  $\bold{F}$
be a non archimedean local field. We assume: 
\ber \label{22}  The characteristic of $\bold{F}$ is different from 2. \eer
Let $\vert.\vert_{\bold{F}} $ be the  normalized absolute value  of $\bold{F}$.

 One considers various algebraic groups defined over $\bold{F}$,  and a sentence
like:  
 \ber \label{corse} " let $A$ be a split torus
 "  will mean '' let 
$A$ be the  group of $\bold{F}$-points   of a torus,  $\underline A$,   defined and split  over $\bold{F}$ ".\eer 
With these  conventions, let $G$ be a connected reductive linear algebraic
group. Let $\tilde{A}_G$ be the maximal split torus of the center of $G$. The change with standard notation will become clear later.
\ste 
Let ${\underline G}$ be the algebraic group defined over $\bold{F}$ whose group of $\bold{F}$-points is $G$. Let $\si$ be a rational involution of ${\underline G}$  defined over
$\bold{F}$.  Let 
$H$ be  the group of $\bold{F}$-points of an open $\bold{F}$-subgroup  of the fixed
point set of $\si$.  We will also denote by  $\sigma$  the restriction of $\si$ to
$G$.\\A split torus  $A$ of $G$ is said $\si$-split if $A$ is contained in  the set of elements of $G$ which are antiinvariant by $\si$. We will denote by $A_G$  the maximal $\si$-split torus of the center of $G$.\\If  $J$ is an algebraic subgroup of $G$ stable by $\si$, one denotes by ${\rm Rat} (J)_\si$ 
the  group of its rational characters defined over $\F$ which are antiinvariant by $\si$.   Let us define:
$${\a}_{G}= {\rm Hom}_{\Z} ({\rm Rat} (G)_\si,\R).$$ The restriction of rational characters 
 from  $G$ to  ${A}_{G}$  induces an  isomorphism: 
\beq \label{iso} {\rm Rat}(G)_\si \otimes _{\Z}\R \simeq {\rm Rat}({A}_{G})\otimes _{\Z} \R .\eeq
Notice that ${\rm Rat} ({A}_G)$ appears as a  generating lattice in the dual space ${\a}'_G$ of ${\a}_G$ and: \beq \label{aotimes} {\a}_G' \simeq {\rm Rat }(G)_\si \otimes _\Z \R. \eeq   
One has the canonical map  ${H}_{G}: G \rightarrow {\a}_{G}$
 which is defined by: 
\beq \label{H} e^{\langle {H}_{G}(x), \chi\rangle}= \vert \chi (x)\vert_{\bold{F}}, \> x\in G,
\chi
\in {\rm Rat} (G)_\si .\eeq The kernel of  
${H}_{G}$, which is denoted by ${G}^1$,   is the intersection of the kernels of the characters of
$G$,  $\vert \chi \vert_{\bold{F}}$, $\chi \in {\rm Rat}
(G)_\si$. One defines
$X(G)_\si= {\rm Hom} (G/{G}^{1}, {\C}^*)$.  It is a subgroup of the group $X(G)$ of unramified characters of $G$. It is precisely the connected component of the neutral element of the group of elements of $X(G)$ which are invariant by $\si$.

   One denotes by  $\tilde{\a}_{G, \bold{F}}$ (resp.,  ${{\a}
}_{G, \bold{F}}$) the image of  $G$ (resp.,  ${A}_{G}$) by  ${H}_{G}$. The group 
$G/{G}^{1}$ is  isomorphic to the lattice 
${\a}_{G,\F}$.\ste 
There is a surjective map: 
\beq  \label{surjection}({\a}'_{G})_{\C}\rightarrow X(G)_\si\rightarrow
1 \eeq   denoted by $ \nu \mapsto \chi_\nu$ which associates to 
$\chi\otimes s$, with $\chi  \in{ \rm Rat}(G)_\si$, $s\in \C$,    the character  $g\mapsto \vert \chi (g)\vert_\bold{F} 
^{s}$ (cf. [W], I.1.(1)).  In other words: 
\beq  \label{cnu}\chi_\nu (g) =   e^{ \langle \nu, {H}_G (g)\rangle}, g \in G, \nu \in({\a}'_{G})_{\C} .\eeq
 The kernel is a lattice and it defines a structure of a
complex algebraic variety on  
$X(G)_\si$ of dimension   
$dim_{\R}{\a}_{G}$.  Moreover $X(G)_\si$ is an abelian  complex Lie group whose Lie algebra is equal to $({\a}'_{G})_{\C}$. \ste
 If  $\chi$ is an element of $ X(G)_\si$, let  $\nu$ be an element of $
{\a}_{G,\C}'$ such that $\chi_\nu= \chi $.  The real part  ${\rm Re}\>\nu \in {\a}'_{G}$ is independent from the choice of $\nu$.  We will denote it  by ${\rm Re} \> \chi$. If  $\chi \in {\rm Hom} (G,
\C^{*}) $ is continuous and antiinvariant by $\si$, the character  of $G$, $\vert \chi \vert$, is an element of  $X(G)_\si$.  One sets ${\rm Re}\> \chi=  { \rm Re}\> \vert \chi\vert $. Similarly,  if  $\chi \in {\rm Hom}({A}_{G},
\C^{*})$ is continuous, the character $\vert \chi \vert$ of ${A}_G$ extends uniquely  to  an element of  $X(G)_\si$  with values in  $\R^{*+}$, that we will denote again by 
  $\vert
\chi \vert$ and one sets  ${\rm Re} \> \chi={  \rm Re}\> \vert \chi
\vert$.\\
A  parabolic subgroup $P$ of $G$   is called a $\sigma$-parabolic subgroup if  $P$ and
$\sigma(P)$ are opposite  parabolic subgroups. Then $M:=P\cap \sigma(P)$ is the  $\sigma$-stable
Levi subgroup of $P$.    If $P$ is such a parabolic subgroup, $P^-$ will denote $\si(P)$. \ber \label{openo} If $P$ is a $\si$-parabolic subgroup of $G$, $PH$ is open in $G$.\eer    The sentence :
 ''Let $P=MU$ be a parabolic subgroup of $G$'' will mean that $U$ is the unipotent radical of $P$ and $M$ is a Levi subgroup of $G$. If moreover $P$ is a $\si$-parabolic subgroup of $G$, $M$ will denote its $\si$-stable Levi subgroup.
 
 If $P=MU$ is a $\si$-parabolic subgroup of $G$,
we  keep the same notations with $M$ instead of $G$.\\ The inclusions
${A}_{G}\subset {A}_{M}\subset M\subset G$  determine a surjective  morphism    ${\a}_{M, \bold{F}}\rightarrow  {\a}_{G, \bold{F}}$ (resp., an injective  morphism, 
 $\tilde{  
{\a}}_{G, \bold{F}}   \rightarrow\tilde{  
{\a}}_{M, \bold{F}}$) which extends uniquely to a surjective linear map between  ${\a}_{M}$ and  ${\a}_{G}$ (resp., injective map,  between
${\a}_{G}$ and 
${\a}_{M}$).
The second map allows to identify  ${\a}_{G}$ with a subspace of 
${\a}_{M}$ and the kernel of the first one, ${\a}^{G}_{M}$,  satisfies: 
\ber $$\label{oplus} {\a}_{M}= {\a}^{G}_{M}\oplus {\a}_{G}.$$
\eer
If an unramified character of $G$ is trivial on $M$, it is trivial on  any maximal compact subgroup of $G$ and on the unipotent radical of $P$, hence on $G$. This allows to identify $X(G)_\si$ to a subgroup of $X(M)_\si$. Then $X(G)_\si $ is the analytic subgroup of $X(M)_\si  $ with Lie algebra ${(\a}'_G)_\C \subset ({\a}'_M)_\C$. This follows easily from  (\ref{cnu}) and (\ref{oplus}).  

Let $P=MU$ be a $\si$-parabolic subgroup of $G$.   Recall that $A_M$ is  the maximal $\si$-split torus of the center of $M$. 
\\Let   $A_P^+,$ (resp., $A_P^{++} $) be the set  of $P$-dominant (resp., strictly
dominant)  elements in $A_M$. More precisely, if  $\Sigma(P)$ is  the
set of roots of $A_M $ in the Lie algebra of $P$, and  $\Delta (P)$ is the
set of simple roots, one has:
$$A_P^+\>(resp., A_P^{++} )= \{a \in  A_M \vert \vert \alpha(a )\vert _{\bold{F}} \geq 1,\>\> 
(resp., > 1) \>\> \alpha \in \Delta (P)\}.$$
Let $A_\vid $ be a maximal $\si$-split torus contained in $M$. Let $\Sigma(U,A_\vid)$  be  the
set of roots of $A_\vid  $ in the Lie algebra of $U$, and  let $\Delta (P, A_\vid)$ be  the
set of simple roots. One defines  for $C>0$ : 
\beqÊ\label{appc} A_\vid ^+(P, C)= \{a \in  A_\vid  \vert \vert \alpha(a)\vert _{\bold{F}}  \geq C , \>\> \alpha
\in \Delta (U, A_\vid)
\}.\eeq 
  Let $A$ be a $\si$-split torus and $g\in G$. We will say that $g$ is $A$-good if and only if  $g^{-1}. A$ is a $\si$-split torus.
Let us prove:
\ber \label{xxsi}If $g$ is $A$-good $\si(g)g^{-1} $ commutes to $A$\eer
It is enough to prove that if $a\in A$,  $(\si(g)g^{-1} ).a= a$. One has
$(\si(g)g^{-1} ).a= \si (g. \si (g^{-1}.a))= \si (g. ( g^{-1} a^{-1})= a$.

For the rest of the article, we  fix $P_\vid= M_\vid U_\vid$ a minimal $\si$-parabolic subgroup of $G$ and let $A_\vid$ be the maximal $\si$-split torus of the center of $M_\vid$. It is a maximal $\si$-split torus of $G$. One denotes by $A_\vid^+$ the set $A_{P_\vid}^+$. A $\si$-parabolic subgroup of $G$ will be said standard (resp., semistandard) if it contains $P_\vid$ (resp., $M_\vid$). 
We choose a maximal split torus $A_0$ which contains $A_\vid$. From [HH], Lemma 1.9, it is $\si$-stable. Let $K_0$ be the stabilizer of a special point of the apartment of the extended building of $G$ associated to $A_0$. 
 \ber \label{WG}From [BD], Lemma 2.4, there exists a finite  set $\mathcal{W}^G_{M_\vid}$ of $A_\vid$-good elements of $G$, such that 
  if $P$ is  any  semi-standard minimal $\si$-parabolic subgroup of $G$, $\mathcal{W}^G_{M_\vid } $ is a set of representatives of the $(P, H)$-double open cosets. We will  assume that $ 1\in\mm{W}^G_{M_\vid}$.  \eer
 For sake of completeness we will recall the definition of $Ê\mm{W}^G_{M_\vid} $. Let $(A_i)_ {i
\in I}$ be a set of representatives of the 
$H$-conjugacy classes  of maximal $\si$-split torus of $G$.  Let us assume that $A_\vid$ belongs to this set. 
The groups  $A_i$ are conjugate under 
$G$ (cf. [HH], Proposition 1.16).
 Let us  choose for each $i$ in $I$, an element 
$x_i$ of  $G$, such that $x_i. A_{\emptyset}  = A_i$  with 
$x_{\emptyset} = 1$. Let $M_i$ be the centralizer of $A_i$ in  $G$. 
If  $L$ is a subgroup of  $G$, one  denotes  by $W_L(A_i)$ the  quotient of  the 
normalizer  in  $L$ of  $A_i$ by its  centralizer. Let us denote by  $W(A_i)$ instead of $W_G(A_i)$.\ste  
 Let  $\mm{W}_{i} $ be a set of representatives in 
$N_G(A_{\emptyset})$ of  $
 W(A_{\emptyset}) / W_{H_i}(A_{\emptyset}) $ where  $H_i = x_i^{-1}. H $. Then   ([HH],  Theorem 3.1)  one can take 
${\cal W}_{M_\vid }^G = \cup_{i \in I} {\cal W}_ix_i^{-1} $.
 
For $ g\in G$ we define $ \dot{g}$ the left coset $Hg $ and we define:
$$\mm{X}^G_{M_\vid} := \{  \dot{ x}\vert x^{-1}  \in \mm{W}^G_{M_\vid} \}.$$
 \section{ The $G$-spaces $X_P$, the constant terms and  the maps $c_{P,Q}$}
\setcounter{equation}{0}
 \subsection{The $G$-spaces $X_P$}
One has (cf. [CD],  Lemma 9.4):
 \ber
 \label{CG94}
 Let $P=MU$ be a $\si$-parabolic subgroup of $G$. The union  of the $(P, H) $ open double cosets in $G$ is equals to $G':= \underline{P}\> \underline{H}  \cap G$. The set $G'$ is also  equal to the set  of $g \in G$ such that $g^{-1}.P$ is a $\si$-parabolic subgroup.
 \eer 
 Let us prove:
\ber  \label{CG} Let $P=MU$ be a $\si$-parabolic subgroup of $G$ and $g \in G$ such that $g. A_M$ is $\si$-split. Then $g.P$ is  a $\si$-parabolic subgroup of $G$. \eer 
 One has  $P=P_\nu$ for some  $\nu \in \a_M'$ in the sense of [CD], (2.14). Then $g.P_\nu= P_{\mu}$ where $\mu$ is the conjugate of $\nu$  by $g$. Our hypothesis implies that $\si(\mu)= - \mu$.  This implies that $g.P_\nu $ is a $\si$-parabolic subgroup    as  $P_\mu$ and $\si(P_\mu)= P_{-\mu}$ are opposite parabolic subgroups.
 
 One easily extends  [CD], Equation (7.1), by replacing $A_\vid$ by $A_M$, the proof being identical:
 \ber  \label{CD71} Let $P=MU$ be a $\si$-parabolic subgroup of $G$. Let $y, y'$ be $A_M$-good elements of $G$ such that $PyH=Py'H$. Then there exist $m\in M, h\in H$ such that $ y'=myh$. 
 \eer 
 
 We define an equivalence relation $\approx_M $ on $\mm{X}^G_{M_\vid}$ by  $x\approx_M x'$ if and only if $Px^{-1}H=Px'^{-1}H$, which by the above equation  is equivalent to $xM=x'M $, as $x^{-1},x'^{-1} $ are $A_\vid $-good. Let $\mm{X}^G_M$  be a set of representatives of the equivalences classes of this  relation. Let  us define $$\mm{W}^G_M: =\{y\in \mm{W}^G_{M_\vid} \vert (y^{-1}) \dot{} \in \mm{X}^G_M\}.$$
 From the above and from (\ref{WG}) one has: \ber  \label{openxp} The set $\mm{X}^G_M$  is a set of representatives of the open $(H, P)$-double cosets in $G$.    \eer
 \begin{lem} \label{XM=}
 Let  $P=MU$ be a semistandard $\si$-parabolic subgroup of $G$. \\ (i) The set of  elements $g$ of $G$  such that $g.A_M$ is $\si$-split  is denoted $X^{Lev}_M\subset G $. It is left  invariant by $H$.  Its  quotient by $H$ on the left is denoted by $X_M\subset H\bb G$.  One has $\mm{X}^G_M\subset X_M$ and: $$X_M = \cup_{x \in \mm{X}^G_M } x M \subset H\bb G , $$  
the union being disjoint.\\
(ii) For each $x \in \mm{X}^G_M $, $xM$ is closed in $X$. 
\\ (iii) We endow $X_M$ with the topology induced by the topology of $X$. Then for each $x \in \mm{X}^G_M$, $x.M$ is open and closed in $X_M$. Moreover   the canonical map 
$(M\cap x^{-1}. H\bb M) \to xM$, $(M\cap x^{-1} .  H)  m\mapsto xm$,  is an homeomorphism.
\\(iv) For  all $x\in X_M$, $xP$ is open in $X$ and $X_MP=X_MU$ is the union of the open orbits of $P$ in $X$. 
\end{lem} 
\dem (i) If $g\in X_M^{Lev}$,  $g.P$ is a $\si$-parabolic subgroup  (cf. (\ref{CG})). From  (\ref{openo})  one has $g^{-1}\in G'$.  One deduces from (\ref{WG}) and the definition of the relation $\approx_M$  that $g^{-1} \in PyH$ for some $y \in \mm{W}^G_M$.  From (\ref{CD71}), one deduces that there exists $m\in M, h\in H$ such that $g^{-1}=myh$. The equality of  (i) follows immediately. From (\ref{CD71}),  if $x,x'$ are distinct elements of $\mm{X}^G_M$,  the sets $HxP$ and $Hx'P$ are disjoint. The disjointness follows.\\
(ii) Changing $H$ into $x^{-1}.H$, one is reduced to prove (ii) when $x$ is equal to ${\dot 1}$. If $(m_n)$ is a sequence in $M$ such that $(\dot{m}_n)$ converges in $X$ to  $l$, then $(\si(m_n)^{-1} m_n)$ 
converges.   The Cartan decomposition for $M\cap H \bb M$ (cf. [BenO]  Theorem 1.1) allows to extract a subsequence of $(m_n)$  denoted again by $(m_n)$ such that $m_n = h_n x a_n \omega_n$, where $(\omega_n)$ converges, $a_n \in A_\vid$, $x^{-1} \in M$ is $A_\vid$-good and $h_n\in M\cap  H$.
Then using (\ref{xxsi}) one has:
$$\si(m_n)^{-1} m_n= \si(\omega_n^{-1})\ \si(x^{-1} )x a_n^2 \omega_n, n\in \N.$$
Hence $(a_n^2)$  is convergent and  $(a_n)$ is  bounded. Extracting again a subsequence we  can assume that $(a_n)$ is convergent. This implies that $(M\cap H ) m_n$ is convergent   in $((M\cap H) \bb M)$ is convergent in $M\cap H \bb M$ and $l$ is element of $\dot{1} M$. This proves (ii).
\\(iii) The fact that $xM$ is closed follows from (ii). As $\mm{X}G_M$ is finite,  (i) implies  that $xM$ is open in $ X_M$. The last assertion follows from [BD], Lemma 3.1  (iii).
\\(iv) From (\ref{CG94}) and (\ref{CG}) and the definition of $X_M$, one sees that $HxP$ is open in $G$. This  achieves to prove the first assertion of  (iv).  The second follows from this and from (\ref{openxp}). \qed
\begin{defi}  Let $P=MU$ be a $\si$-parabolic subgroup of $G$. 
 Then $X_M$ is a $P^-$-space with the given action of $M$ and with the trivial action of $U^-$. 
 We define: 
 $$X_P= X_M \times_{P^-} G.$$
Then $X_M$ identifies to a subset of $X_P$.  If $x\in X_M$, its image in $X_P$ will be denoted by ${x}_P$.
\end{defi}
If $x, x' \in X_M$ the notation 
$x\approx _Mx'$ will mean  that $x,x'$ are in the same $M$-orbit in $X_M$. The following assertion follows from the definition of $X_P$. 
\ber  \label{xpg} Let $x, x' \in X_M$. 
The following conditions are equivalent:\\
(i) $x_PG=x'_P G$.\\
(ii) $xM=x'M$ in other words $x\approx_M x'$.
\eer
We define $H_P:= U^-(M\cap H)$. If $y \in G$, let us denote by $\si_y$ the rational involution of $G$ defined by:
$$\si_y(g)= y^{-1}  \si   (ygy^{-1}) y,$$ whose fixed point set is equal to $y^{-1} .H$  Moreover $\si_y$ depends only on $\dot{y}$. 
\ber \label{stab}Let $x\in X_M \subset H\bb G$. The stabilizer of ${x}_P $ in $G$ is equal to $ (x^{-1}.H)_P:= U^-(M\cap x^{-1}.H). $\eer

\begin{defi} \label{aaction}  Let $a\in A_M$. From (\ref{xpg})  any element  $y\in X_P$ is  of the form $y= x_P g$ for  a unique element  $x \in \mm{X}^G_M $ and some element $ g \in G$, which is defined up to the left action of $U^- (M\cap x^{-1}.H)$ .  We   see easily from (\ref{stab}) that  $a. y:= x_P a g $ is well defined. It defines a left  action of $A_M$ on 
$X_P$ which commutes to the right $G$-action. \end{defi} 
From the equality in Lemma \ref{XM=}, one deduces the following equality:
\begin{lem}  \label{disun} 
(i) One has:  \beq X_P = \cup_{x\in \mm{X}^G_M } {x}_P G,  \eeq 
the union being disjoint.\\
(ii) For $x\in \mm{X}^G_M $, $x_PP$ is the unique open orbit in $x_PG$.\\
(iii) Let $(X_M)_P$ the image of $X_M$ in $X_P$ or equivalently the set $\{x_P\vert x \in X_M\}$.  The union of the open $P$-orbits in $X_P$ is equal to $(X_M)_P P= (X_M)_P U$ and the map  from $X_MP= X_M U$ to $(X_M)_P P$ defined by 
$xu \mapsto x_P u$ is a bijective $P$-equivariant map. 
\end{lem}
\dem
One deduces (i) from the equality in Lemma \ref{XM=}.
\\(ii)  It follows from (\ref{stab}) that  $x_PG$ is isomorphic to $ U^- (M\cap x^{-1}. H) \bb G$. Then (ii) follows from the fact that there is a unique open  $(U^-, P)$-double coset in $G$.
\\ The first part of (iii) is clear.  It follows from  (\ref{stab})   that  the map $X_M \times U\to (X_M)_P P$, $(x,u) \mapsto (x_Pu)$ is bijective. One checks easily that it is $P$-equivariant. \qed 
Let  $P$ be  a standard $\si$-parabolic subgroup of $G$. Let us prove:  
\beq \label{aina}\{ a \in A_\vid\vert \vert \aa (a) \vert _\F \geq C, \aa \in \Delta(P_\vid \cap M)\}  = A^+_\vid(P_\vid, C) A_M .\eeq  
The right  hand side is clearly included in the left  hand side of the equality to prove. Let $ a $ be an element of the left  hand side. Let $b \in A_M$ be strictly $P$-dominant. Then for large $n\in \N$, one has $ab^n \in A^+_\vid(P_\vid, C) $. Our claim follows. 
\begin{prop} \label{omega}
There exists a compact subset $\Omega$ of $G$ such that for all $\si$-parabolic subgroup $P$  of $G$ containing $M_\vid$, one has:
$$X_P= \cup_{x \in \mm{X}^G_{M_\vid}} x_P A_\vid^+A_M\Omega. $$
\end{prop} 
\dem The claim is true for $P=G$ from the Cartan decomposition for symmetric space  (cf. [BenO]  Theorem 1.1).
 In general one has $G=P^-K_0$ hence $$X_P= X_M P^-K_0= X_M K_0= \cup_{x\in \mm{X}^G_M} x_P M K_0.$$
 The $M$-space $xM\subset H \bb G$ is a symmetric space for $M$ for the  involution $\si_x$ restricted to $M$. As   $x$ is $A_\vid $ good,  
$P_\vid\cap M$ is a $\si_x$-parabolic subgroup of $M$ (cf. [CD] Lemma 2.2).  From the Cartan decomposition for this  symmetric space, it is enough to prove the following lemma.
\begin{lem}
The  open orbits of $P_\vid\cap M$ in $xM$ are the orbits $y(P_\vid\cap M)$, where $y$ describes the set of elements in $\mm{X}^G_{M_\vid} $ such that $y\approx_M x$.
\end{lem}
By conjugating on the left by $x^{-1}$ and changing $H$ into $x^{-1}. H$ one is reduced to prove the lemma for $x=\dot{1}$.  Any open $(P_\vid \cap M)$-orbit in $(M \cap H) \bb M$ is of the form 
$  (M\cap H) z (P_\vid\cap M)$  where $z^{-1}  $ is $A_\vid$-good and element of $M$ (cf. (\ref{WG})). 
 As $HP$ is open, the  product map $H\times P \to H P$ is open  (cf. [BD], Lemma 3.1 (iii)). Hence, as $H z P_\vid =H ((H\cap M) z (P_\vid\cap M))U$,  one sees that $H z P_\vid $ is open.  Then  (\ref{WG})  implies the existence of  an element  $y$ of $ \mm{X}^G_{M_\vid} $  such that: \beq \label{yz} HzP_\vid= HyP_\vidÊ.\eeq  As $z\in M$,  $z  $ is $A_M$-good. As $y$ is also $A_M$-good, it follows from (\ref{CD71}) that  $z= hym' $ for some $m'\in M$, $ h \in H$ and one has  $y\approx_M z$. Let us prove: $$(Hz P_\vid )\cap H M  = Hz(P_\vid \cap M).$$ Let $p\in P_\vid $ and let us write $p= p'u$ with $p' \in P_\vid \cap M$ and $u \in U$. Let us show that $zp' u\in HM $ if and only if $ u=1$. Let $m':= zp'\in M$.  If  $z p' u\in HM$, there exist $h\in H$, $m\in M$ such that   $m'u= hm $.  Then one has
$$h = m'm^{-1} (m.u)$$ Hence both sides of the equality are elements of $H\cap P= H\cap M$. It follows that $u=1$. Our claim follows. 

 As $z\in M$,  $\dot{z}\approx_M \dot{1}$. 
Taking into account $y \approx_M z $,  one has $y \approx_M \dot{1}$   and one shows  similarly that: \beq  \label{yinter} (Hy P_\vid) \cap H M = Hy(P_\vid \cap M).\eeq   
From this  and (\ref{yz}) one sees that: $$Hy(P_\vid \cap M)= Hz(P_\vid \cap M).$$
This shows that any open $P_\vid \cap M$-orbit   in $\dot{1} M$ has the required form.  

Reciprocally from (\ref{yinter}) one sees that for all $y\in \mm{X}^G_{M_\vid}$ such that $y\approx_M \dot{1}$, $  y(P_\vid \cap M) $ is open in $X_M$ as it is equal to the intersection of  an open set of $X$ with the open subset $\dot{1} M$ of $X_M$ (cf. Lemma \ref{XM=} (iii)). This proves the Lemma. \qed 
\begin{rem}
1) There is a minor change with [SV]. Here we are interested to $X= H \bb G$ but Sakellaridis and Venkatesh study the bigger space $(\underline{H} \bb \underline{G})(\F)$. The space $X$ appears as one of the finitely many $G$-orbits in $\underline{X}(\F)$ and every $G$-orbit in  $\underline{X}(\F)$ is of the same type than $X$.
\\ 2) If $P=MU$ is a   standard $\si$-parabolic subgroup of $G$, we define  $\Theta_P$   as  the set of  simple $A_\vid$-roots  in the Lie algebra of $M$  which are simple for $P_\vid$. Notice that $\Theta_{P_\vid} = \vid$. We could define also $A_{\Theta_P}=A_P$. Then $A_{\Theta_P}$ plays here the role $A_{X, \Theta_P}$ in [SV]. 
\end{rem}
  \subsection{Constant term}
  
  Let $J$  be a totally discontinuous group  acting continuously on a totally disconnected topological space $Y$. We will say that the action is smooth if the stabilizer of any element of $Y$ is open and  we will denote by $C^\infty(Y)$ the space of functions which are fixed by the right action of some compact open subgroup of $G$. 
  
  Let us recall (cf.  [D], Proposition 3.14) the following result. 
   \ber \label{ctepp}
   Let $P=MU$ be a $ \si$-parabolic subgroup of $G$. Let $(\pi,V)$ be a smooth $G$-submodule of $C^\infty (H\bb G)  $. The map $fÊ\to f_P$ is
the unique morphism of $P$-modules from $V $ to the space  $ C^\infty((M \cap H)Ê\bb M)$  endowed with the  right 
action of $M$ tensored by $\delta_P^{1/2}$ and the trivial action of $U$,  such that:\\
 For all  compact
 open subgroup, $J$, of $G$ there exists  $C>0$ , such that for all  $f \in V^J $:
 $$  f(a)=  \delta_P^{1/2} (a )f_P(a), a\in A_M(P^-, C)= A_M \cap  A_\vid ^+(P^-, C),  $$
 where $A_\vid ^+(P^-, C)$ has been defined in (\ref{appc}).\\
 We have a similar statement by replacing the preceding equality by 
 $$  f(a)=  \delta_P^{1/2} (a )f_P(a),  a \in A_\vid ^+(P^-, C) $$
 \eer 
We have slightly modified the statement of l.c.   by replacing $A_0 $ by $A_M $ and $A_\vid$  but unicity still holds due to  [D] Equation (3.8). It is useful to introduce:  \beq \label{tf} \tilde{f} _P= \delta_P^{1/2}f_P.\eeq 
\ber \label{rhoafp} Let us assume that $V$ is of finite length. Let $(\delta,E)$ be the unormalized Jacquet module of $V$. Then there exists a finite family of complex characters $\chi_1,\dots, \chi_r$ of $A_M$ such that
 $$(\delta(a)-\chi_1(a))\dots (\delta (a)-\chi_r(a))=0, a\in A_M$$
 From the interwining properties of the map $f\mapsto f_P$, one deduces $$ (\rho(a)-Ê\chi_1(a))\dots (\rho(a)-Ê\chi_r(a))\tilde{f}_P= 0, a \in A_M.$$\eer
\begin{theo} \label{theoct} Let $P=MU\subset Q=LV $ be two  standard $\si$-parabolic subgroups of $G$. If $C\geq 0$, let $A^+_\vid(P,Q, C)$ be the set of $a\in A_\vid^+ $   such that $ \vert \aa(a) \vert_\F \geq C $  for all roots $ \aa$  of $ A_\vid $  in  the   Lie algebra of $U\cap L$. 
\\ (i) There exists a unique $G$-equivariant  map $c_{P,Q}$  from $C^{\infty}(X_Q) $ to $C^\infty(X_P)$ satisfying the following property:
\\For all compact open subgroups $J$ of $G$, there exists $C>0$ such that  for all $f \in C^\infty(X_Q) $ which is right $J$-invariant, one has:
\ber \label{cpqdef}  
$$(c_{P,Q}f) ( x_P a  )= f( x _Q a  ),a\in A_\vid^+(Q, P,C), \>  x  \in \mm{X}^G_{M_\vid}.$$
\eer
The map  does not depend on the choice of $\mm{W}^G_{M_\vid}$.
\\(ii)  Let $R$ be an other standard $\si$-parabolic subgroup of $G$ such that $Q\subset R$. Then one has:
$$c_{P, R}= c_{P,Q} \circ c_{Q,R}.$$
\\ (iii) Let $\mm{V}$ be a smooth $G$-submodule of finite length of $C^\infty(X_Q)$. Then there exists a finite family of complex characters $\chi_1,\dots, \chi_r$ such that for all $f \in \mm{V}$:
$$((\l(a)-Ê\chi_1(a))\dots (\l(a)-Ê\chi_r(a))c_{P,Q}f)(x_Pg)= 0, x\in \mm{X}^G_M, g \in G, a \in A_M.$$
\end{theo}
For the proof we will need two lemmas.
\begin{lem} \label{cpqx}
Let $x\in X_M$. \\
(i)If $f \in C^\infty(X_Q) $ and $g\in G$, let  $f_{x_Q,g}$ be the map  $l \mapsto f(x_Q  l g) $ viewed as a map on $ (x^{-1}.H)\cap L\bb L$. We define a function  
 $f_{x_Q, P^-\cap L}$ on $G$ by $g \mapsto (f_{x_Q,g})^{\tilde{}}_ {P^-\cap L}(  1)$,  where  we use the notation (\ref{tf}). 
 It  is left invariant by $(x^{-1}.H)_P$ and it is right $J$-invariant if $f$ is right $J$-invariant.\\
(ii) The point (i)  allows to define a map $c_{P,Q,x}: C^\infty (X_Q)\to C^{\infty} (x_P G)$ by
$$(c_{P,Q,x}f)(x_P g)= f_{x_Q, P^-\cap L}(g).$$
It intertwines   the right regular representations of $G$ on $C^\infty (X_Q)$ and  $C^{\infty} (x_P G)$. \\
(iii) 
One has 
$$(c_{P,Q ,x}f) (x_P m g) := (f_{x_Q, g})^{\tilde{}}_ {P^-\cap L}(  m), m \in M.$$
(iv) For all compact open subgroup $J$ of $G$, there exists $C>0$ such that for all $x  \in \mm{X}^G_{M_\vid}$, for all $f \in C^{\infty}(X_Q) $ which is right $J$-invariant, one has:
  $$(c_{P,Q ,x}f)( x_P a  )= f (x _Q a ),a\in A_\vid ^+(P,Q,C), \>  x  \in \mm{X}^G_{M_\vid}.$$
(v)  We have unicity of the $G$-maps satisfying the condition above on the sets $A_M\cap A^+_\vid(P,Q, C) $.
\end{lem}
\dem
 (i)  Due to  the intertwining properties of the constant term map (cf. (\ref{ctepp})) the map $\varphi\mapsto \tilde{\varphi}_{P^-\cap L}$ intertwines the right  regular representations of $P^-\cap L$ on   
$C^\infty((x^{-1}.H)\cap L)\bb L)$ and on  $ C^\infty ((x^{-1}.H) \cap M)\bb M) $, where $U^-\cap L$ acts trivially on the latter space. Also one remarks that $f_{x_Q, vg}=f_{x_Q, g}$ for $g\in G, v\in V$. 
Altogether this shows (i) and that the map $c_{P,Q,x} $ is well defined. The map $c_{P,Q,x} $ intertwines the right regular representations of $G$ as the equality  $(c_{P,Q,x}(\rho(g)f))(x_P g')= (c_{P,Q,x}f)(x_P g'g)$ follows from the definitions.  This achieves to prove (ii). \\
 (iii) By (ii), it is enough to prove this for $g=1$. The intertwining  properties of  the map $\varphi\mapsto \tilde{\varphi}_{P^-\cap L}$  described above allows to prove (iii). \\
 (iv) By (iii) and from the second equality of (\ref{ctepp}) for $P^-\cap L$, one deduces (iv).\\
 (v) As $c_{P,Q,x}$ is a $G$-map,  (v) follows from the second characterization in  (\ref{ctepp}) of the constant term.\qed  
\begin{lem} \label{cpqxy}
Let $ x, y \in X_M$. If $x\approx_My$, one has $c_{P,Q,x}=c_{P,Q,y}$.
\end{lem}
From  Lemma \ref{cpqx} (v), it is enough to prove the following assertion.  
\ber Let $J$ be a compact open subgroup  of $G$. There exists $C>0$ such that for all $f \in C^\infty (X_Q)$  which is $J$-invariant
$$ \label{cpqy} (c_{P,Q,x} f)(y_P a )=f(y_Q a), a\in A_M\cap A_\vid^+(P,Q,C).$$
\eer   
Let $ m\in M$ such that  $y=xm$. Then $y_P= x_P m, y_Q= x_Q m $. 
By the interwinining properties of $c_{P,Q,x}$ and the commutation of $a\in A_M $ with  $m$, 
one has \beq \label{cpqrho} (c_{P,Q,x} f)(y_P a ) = c_{P,Q,x} (\rho({m}) f)(x_P a ), a \in A_M.\eeq 
 One remarks that $\rho({m})f $ is fixed by $m.J$.  Hence as $x$ satisfies  Lemma \ref{cpqx} (iv) , there exists $C>0$ such that for all $f \in C^\infty (X_Q)$ right invariant by $J$:
$$ c_{P,Q,x} (\rho({m})f)(x_P a ) = (\rho({m}) f)(x_Q a), a\in A_M\cap A_\vid^+(P,Q,C).$$ 
As $(\rho(m)f)(x_Q a)= f(y_Q a)$, together with (\ref{cpqrho}) this proves (\ref{cpqy}) and the lemma. \qed 
\\{\em  Proof of Theorem 1}\\
(i) We define $c_{P,Q} (f)$ for $f \in C^\infty (X_Q)$ by: 
$$(c_{P,Q}f)(x_Q g):= (c_{P,Q, x} f )(x_Q g), x \in \mm{X}^G_M$$
From Lemma \ref{cpqx} (iv) and (v), one sees that this is well defined and that it has the required properties including unicity.
Also from Lemma  \ref{cpqxy}, $c_{P,Q}$ does not depend on the choice of   $\mm{X}^G_M$ in $\mm{X}^G_{M_\vid} $. Also,  as changing our choice of $\mm{X}^G_{M_\vid} $ involves only right multiplcation by elements of $M_\vid$, one sees that $c_{P,Q}$ even does not depend of the  choice of  $\mm{X}^G_{M_\vid}$. \\(ii)   follows  easily from the unicity statement in (i).\\(iii) We use the notation of Lemma \ref{cpqx} (i). The map $f\mapsto f_{x_Q,1}$ is a  $Q^-$-map from $\mm{V}$  to a $Q^-$-submodule of $C^\infty((x^{-1}.H)\cap L\bb L)$ endowed with the right action of $L$ and the trivial action of $V$. This submodule is  a quotient of the unormalized Jacquet module of $\mm{V}$ for $Q^-$. Hence it is an $L$-module of finite length. Then (iii)   follows from the definition of $c_{P,Q}$ above and from (\ref{rhoafp}) applied to $L$ instead of $G$. 
\qed
\section{Neighborhoods at infinity of $X_Q$ and mappings $exp_{X_P, X_Q}$} 
\setcounter{equation}{0}
 \subsection{Choice of measures} \label{mea}
We fix on $G$ (resp., $H$, resp.,  the unipotent radical of a semistandard $\si$-parabolic $P=MU$   of $G$) the Haar measure such that its  intersections with $K_0$ is of volume 1. 
From this we deduce  a measure on $H\bb G$.
We choose the Haar measure on $M$   such that:  \beq \label{intumu}\int_G f(g)dg= \int_{U\times M \times U^-} f(umu^-) \delta_P(m)^{-1} dudmdu^-, f \in C_c^\infty(G).\eeq 
Also there exists a constant $\gamma(P)$ such that: 
\beq \label{umk} \int_G f(g)dg =\gamma(P)  \int_{U^-\times M\times K_0} f(u^-mk) du^- dm dk . \eeq   
The set $X_MU$ is an open subset of $H\bb G$ (cf. Lemma \ref{XM=} (iv)) which is right invariant by $P$. Hence the measure on  $H \bb G$ induces a right $P$-invariant measure on $X_MU$. But  the map $X_M \times U \to X_MU$, $(x,u)\mapsto xu$ is a homeomorphism. As the Haar measure on $U$ has been fixed, there is a canonical measure $m_{X_M}$ on $X_M$ such that:
 \beq \label{dmx}  \int_{X_M U} f( y) dy= \int_{X_M \times U } f(xu) dm_{X_M}(x) du, f\in C_c(X_P). \eeq 
 One checks easily that this measure satisfies:
 \beq \label{intdx} \int_{X_M} f(xm)dm_{X_M}(x)= \delta_P(m)^{-1}\int_{X_M}  f(x)dm_{X_M}(x), m\in M \eeq
  Let $x\in X_M$. As $U^-P$ is open in $G$, $x_P P$ is an open set in $X_P$ which depends only on $xM$.   By looking to the stabilizer of $x$ and $x_P$ one sees that the map $xp\mapsto x_Pp$ is a well defined continuous bijection between $x P$ and $x_PP$ which depends only on $xM$ hence  on $x_PP$. Thus, our choice of $P$-invariant  measure on ${x}P$ induces  and ''by transport de structure'' a $P$-invariant measure  on  $ x_PP$. We fix on $x_P G$ the $G$-invariant measure which agrees with this measure on ${x}_P P$.  Hence we have a right invariant measure by $G$ on $X_P$.
 We want to deduce from $m_{X_M}$  an $M$-invariant measure on $X_M$. This will depend on our choice of $\mm{X}^G_M$.
If $x \in \mm{X}^G_M$, the map $(M\cap x^{-1}. H)\bb M\to xM$, $(M\cap x^{-1}. H)m \mapsto xm$ is a homeomorphism (cf. e.g. [BD] Lemma 3.1 (iii)). The measure on $X_M$ determines a measure on $(M\cap x^{-1}. H)\bb M$.  Let us show:
\ber  \label{dptriv} The function $\delta_P$ is trivial on $M \cap  x^{-1}.H$. \eer The group  $P$ is a  $\si_x$-parabolic subgroup of $G$  (cf. [CD], Lemma 2.2 (iii) where one has  to change $x$ in $x^{-1}$). This implies that  $\delta_P$ is antiinvariant by $\si_x$  and hence trivial on the fixed points of $\si_x$. 
This proves our claim. 
This determines  `par transport de structure'' a function denoted $\delta_{P, x} $ on $xM$. Multiplying the  restriction to $xM$ of the canonical quasiinvariant measure $m_{X_M}$ by $\delta_{P,x}$ one gets an $M$-invariant measure on $xM$ and on $(M \cap x^{-1}.H) \bb M$. 
Hence one has: 
\ber \label{dxm} Our choice of $\mm{X}^G_M$ determines an $M$-invariant measure  on $X_M$. \eer
It allows to identify $ C^\infty(X_M)$ to  a subspace of  the  dual of $ C_c^\infty (X_M)$ (we will see later that this subspace of the dual is the full smooth dual, cf. after (\ref{sdual})).

One deduces also a measure on $x^{-1}.H$ by conjugacy.  Together with our choice of measure on $M$ and on $(M \cap x^{-1} .H) \bb M$, this   determines a measure on $(M \cap x^{-1}. H) \bb x^{-1}. H$.
 
 We introduce a unitary  action  $\mm{L}$  of $A_M$ (cf. (\ref{intdx}) for unitarity)  on the space $L^2(X_P)$ called normalized action:
 \beq  \label{mml} \mm{L}_a f(x) =  \delta_P^{1/2} (a) f(ax), x\in X_P, Ê\eeq
 where $ax$ is the left action of $a\in A_M$ on $x\in X_P$ of  Defintion \ref{aaction}.

\subsection{Compact open subgroups with a $\si$-factorization}
First we  give   a definition. \ber \label{sifac}  A compact open subgroup $J$ of $G$ is said to have a $\si$-factorization (resp. strong $\si$-factorization) for standard $\si$-parabolic subgroups of $G$ if it satisfies the following conditions:\\ 
(i)  For every standard  $\si$-parabolic subgroup $P=MU$  of $G$ the product map $ J_{U^-}\times J_M \times J_U \to J$ is bijective, where  $J_{U^-} = J \cap U^-$, $J_M=
J\cap M$,  $J_U= J \cap U. $\\
(ii)  Let $A\subset A_\vid$ be the maximal $\si$-split torus of the center of $M$ and let $A^-$ (resp. $A_\vid^-$) be the set of its $P$-(resp. $P_\vid$)-antidominant elements.  For all $a$ belonging to  $A^- $ (resp. $A_\vid^-$ for the strong $\si$-factorization) one has $$aJ_U a^{-1}  \subset J_{U},  a^{-1}J_{U^-} a \subset J_{U^-}.$$
(iii) One has $J= J _H J_P$, where $J_H= J \cap H, J_P= J \cap P$. \\(iv) For every  $\si$-parabolic subgroup $P=MU$ of $G$ which contains $P_\vid$,   $J\cap M$ enjoys the same properties that $J $ for $M$ and  $P_\vid\cap M$. \eer
From [CD]  Prop 2.3, there exist arbitrary small compact open subgroups of $G$ with a $\si$-factorization.
We will need the following lemma later.
\begin{lem} \label{strongsi}
 There exists a basis of neighborhood of the identity  in $ G$, $(J'_n)_{n \in \N}$,  made of  a decreasing  sequence of compact open subgroups of $G$ with a  strong $\si$-factorization and such that for all $n\in \N$, $J'_n$ is a normal subgroup of $J'_0$.
\end{lem} \dem
We keep the notation of  [CD]  Prop 2.3, Then,  as the basis of $\un{u}_\vid$ and $\un{u}^-_\vid$ is made of weight vectors $\un{a_\vid}$,  one has: $$\LL\underline{g}=\LL \un{u} \oplus\LL \un{m} \oplus\LL \un{ u^-},$$ where $\LL\un{u}= \LL\un{g} \cap \un{u}$, $\LL\un{m}= \LL\un{g} \cap \un{m}$, $\LL\un{u}^-= \LL\un{g} \cap \un{u}^-$ and $\LL \un{u}$ (resp., $\LL \un{u}^-$ ) is stable by the adjoint action of $A_\vid^-$ (resp., $A_\vid^+$). Then one shows  as in  the proof of  [CD]  Proposition  2.3, where only (ii) has to be modified,  that there exists a basis of neighborhoods $(J_n)_{n \in \N}$ of the identity  in $ G$   made of a decreasing sequence of compact open subgroups of $G$ with a  strong $\si$-factorization. 

As $\LL\un{g} $ is  compat and open in $\un{g}$, there exists $n_0\in \N$ such that  the adjoint action of $J_{n_0}$ preserves $\LL\un{g}$. Hence by l.c. Lemma 10.1 (iii), there exists $N\in \N$ such that  for all $n $ greater than $\N$, $J_{n_0}$ normalizes $J_n$. The sequence $(J'_n)$ defined by  $J'_n= J_{N+n}$ has the required properties.\qed 
\subsection{Statement of Theorem \ref{theoexp}}
\begin{theo}  \label{theoexp} Let $P=MU\subset Q=LV$ two standard $\si$-parabolic subgroups of $G$. Let    $K$ be a    compact open subgroup of $G$  having a strong  $\si$-factorization. Let $\Omega$ be as in Proposition \ref{omega}.   We may and will assume that $K\subset \Omega $ and that  $\Omega$ is  biinvariant by $K$.   Let  $J$ be  a compact open subgroup of $G$ such that for all $\omega$ in $\Omega$, $x^{-1} \in \mm{W}^G_{M_\vid}$ , $(x\omega).  J \subset K$. 
 
 We define    for $C>0$ and $x \in \mm{X}^G_{M} $: $$N_{X_Q} (x, P,C): =\cup_{ y \in \mm{X}^G_{M_\vid} , y \approx_M x} y_Q A_{\vid} ^+ (P, Q,C) \Omega.$$ Then there exists  $C>0$  such that:\\ (i) The union  $$N_{X_Q  } (P,C) := \cup_{x\in  \mm{X}^G_M } N_{X_Q}(x,P, C)$$ is disjoint. \\
     (ii) The subset $N_{X_Q } (P,C) $ of $X_Q$ is right $J$-invariant. We view $N_{X_Q , J } (P,C): =N_{X_Q }(P, C) /J  $   as a subset of $  X_Q  /J$. The map  $N_{X_Q,J}(P,C)  \to X_P/J$  which associates $x_P a  \omega J$   to ${x_Q } a \omega J$ with $x \in \mm{X}^G_{M_\vid} $, $a \in A_\vid ^+(P,Q, C)$, $\omega \in\Omega$  is well defined on  $N_{X_Q , J } (P,C)$.  It is denoted $exp_{X_P, X_Q, J} $.The image  by $exp_{X_P, X_Q, J} $ of $N_{X_Q , J } (P,C)$  is equal to  $N_{X_P , J } (P,C)$\\
     (iii) The map $exp_{X_P, X_Q, J} $ is injective on  $N_{X_Q , J } (P,C)$. \\
 (iv) As a map from  a set of $J$-invariant  subsets  of   $X_Q$ to  a set of  $J$-invariant  subsets of $X_P$,  $exp_{X_P,X_Q,J}$  preserves volumes.
 
\end{theo}
From the definitions, one sees:
\begin{cor} \label{cortheoexp}
If $a \in A_L$ is $Q$-dominant  and $z\in N_{Q,J}(P,C)$, one sees from the definitions that $az\in N_{X_Q,J}((P,C)$ and that:
$$exp_{X_P, X_Q,J}(az)= aexp_{X_P, X,J} (z), z \in  N_{X_Q , J } (P,C)$$
\end{cor}
\setcounter{cor}{0}
{\em  First reduction for the proof of Theorem \ref{theoexp}}\\
We will reduce the proof of the theorem to the case where $Q= G$.  The proof when $Q=G$ will be done in section \ref{endtheo}.  Let us assume that the theorem has been proved for $Q=G$. Let us prove it for arbitrary $Q$.

We will define $exp_{X_P, X_Q,J}$ and prove part (ii) of Theorem \ref{theoexp}.  We define $N'_{X_Q ,J } (P,C) = exp_{X_Q, X, J}(N_{X ,J }(P,C))$ which is well defined for $C$ large. 
Then, from (\ref{aina}),  the definition of the left $A_L$-action (cf. Definition \ref{aaction})  and the definition of $exp_{X_Q,X,J}$ one has: $$N_{X_Q , J } (P,C)= A_L N'_{X, J } (P,C).$$
Let $y \in N_{X_Q , J } (P,C)$. By the above equality and the definition of $ N'_{X_Q , J } (P,C)$, there exist  $a \in A_L$ and $z \in N_{X , J }(P,C)$ such that $y = a exp_{X_Q, X,J} (z)$. Let us  prove that $aexp_{X_P, X,J} (z)$ does not depend on the choice of $a$ and $z$ as above.

Let us assume that there exists $a'\in A_L$ and $z' \in N_{  X, J }(P,C)$ with $y= a' exp_{X_Q,X, J} (z')$. By choosing $b\in A_L$ sufficiently $Q$-dominant we can assume that $ba, ba' $ are $Q$-dominant. 
As $z \in N_{X,J}(P,C)$ one may write  $z= x a_z \omega J$ for some $x \in \mm{X}^G_{M_\vid}$, $a_z \in A_\vid^+(P,C)$, $\omega\in \Omega$. By abuse of notation, as it may depends on this writing, one defines
$baz:= x ba a_z \omega J$. One defines similarly $b'a'z'.$
Then $baz, ba'z'\in N_{X, J }(P,C)$. 
From our hypothesis one has: $$ba exp_{X_Q, X,J} (x)= ba'exp_{X_Q, X,J} (x').$$
From Corollary \ref{cortheoexp} of Theorem \ref{theoexp} for $Q=G$,  one has:
$$ba exp_{X_Q, X,J}( z)= exp_{X_Q, X,J}  (baz), ba' exp_{X_Q, X,J}( z)= exp_{X_Q, X,J}  (ba'z').$$
From the injectivity in (iii) for $Q=G$, one deduces: $$baz= ba'z' .$$ One sees from the definition of $exp_{X_P,X,J}$ in (ii) that:
$$exp_{X_P, X,J}  (baz)= ba exp_{X_P, X,J}( z), exp_{X_P, X,J}  (ba'z')= ba' exp_{X_P, X,J} (z').$$
As $baz=ba'z'$, one deduces from this the equality:
$$a exp_{X_P, X,J}( z) = a' exp_{X_P, X,J} (z').$$
This proves our claim and it 
 allows to define $$exp_{X_P, X_Q,J}(y):= aexp_{X_P, X,J} (z).$$
Let $y= x_Qa\omega J\in N_{X_Q,J}(P,C)$ with $x \in \mm{X}^G_{M_\vid}$, $a \in A_\vid^+(P,Q,C)$.  By 
choosing $b' \in A_L$ sufficiently $Q$-dominant,  one has $a':= b'a \in A_\vid^+(P,C)$. Let $b= b'^{-1}$   and $y'=  x_Qa'\omega J$
One has $y= by'$ and $y'= exp_{X_Q,X,J} (xa'\omega J)$. Our definition of $exp_{X_P, X_Q, J}$ shows that:
$$exp_{X_P, X_Q, J}(x_Qa\omega J)=b x_P a'\omega .J= x_P a \omega J.$$

This achieves to prove  that  $exp_{X_P, X_Q, J} $ is defined by the formula given in the theorem.  This implies that the  the image of $N_{X_Q , J } (P,C)$ is clearly $N_{X_P , J } (P,C)$.  This achieves the proof of Theorem \ref{theoexp} (ii)  and Corollary  \ref{cortheoexp} follows.  \\ 
(iii) Let $y, y' \in N_{X_Q,J}(P,C)$ with  $exp_{X_P, X_Q, J}(y)= exp_{X_P, X_Q, J}(y')$. One wants to prove that $y=y'$. By multiplying $y$ and $y'$ by a sufficiently $Q$-dominant element of $A_L$, 
one may assume that $y, y'\in N'_{X_Q,J} (P,C)$.  Then $y = exp_{X_Q, X,J} (z), y' = exp_{X_Q, X,J}(z')$ with $z, z' \in N_{X, J}(P, C)$.  From our definition of $exp_{X_P, X_Q,J}$, one deduces the equality:
$$exp_{X_P, X,J}(z)= exp_{X_P, X,J}(z').$$
From the injectivity of $exp_{X_P, X,J}$ one sees that $z=z'$, hence $y=y'$.
This achieves to prove (iii).
\\(iv) One has the equality $$vol_{X_P}  (axJ)= \delta_P(a) vol_{X_P}  (x J), x\in X_P, a \in A_P.$$
Using this equality for $P$ and $Q$, using Theorem \ref{theoexp} for $Q=G$ and $P$ successively equal to $P$ and $Q$, and our definition of $exp_{X_P,X_Q,J}$ one deduces (iv) for all $Q$.

It remains to prove (i). One has  $y_PG=x_PG $  if  and only if $x\in \mm{X}^G_M $ and $y \in X_M$ is  such that $x\approx_M y$ (cf. Lemma \ref{XM=}). From the ''if part'' and the definition above of $exp_{X_P,X_Q,J}$, the image of $N_{X_Q,}(x, P, C)$ by $exp_{X_Q,X_P, J} $  is contained in  $x_PG$. Then the ''only if part ''   implies (i). \qed 

The following proposition is an easy consequence  from the definition in part  (ii) of the Theorem  above. 
\begin{prop}
With the notation of Theorem, \ref{theoexp}, one has
$$exp_{X_P, X_Q,J} ( exp_{X_Q,X,J } (xJ))  = exp_{X_P,X,J }(xJ), x \in N_{X,J} (P,C)$$
\end{prop}
The following assertion is an immediate corollary of the Cartan decomposition for $X_Q$.
 \ber \label{unionn}
 Let $C>0$. The complementary set in $X_Q$ of the union of $N_{X_Q  } (P,C) $ when $P$ describes the maximal standard $\si$-parabolic  is a compact set modulo the action of $A_L$. 
 \eer
 
\section{Eisenstein integrals and some results of Nathalie Lagier}
\setcounter{equation}{0}
\subsection{Eisenstein integrals}
Let $P=MU$ be a semi standard $\si$-parabolic subgroup of $G$. 
Let $(\delta, E)$ be a unitary irreducible smooth representation of $M$. Let $\chi \in X(M)_\si$ and let $\delta_\chi= \delta\otimes \chi$ and let us  denote  by $E_\chi $ the space of $\delta_\chi$. Let $(i^G_P \delta_\chi ,  i^G_PE_\chi) $ be the  normalized parabolically  induced representation.  
 
\ber \label{duali} The intertwining linear map from   $i^G_P\cc{\delta}$ to $ (i^G_P \delta)^{\cc{} }$  which associates to $\cc{v} \in i^G_P\cc{\delta} $ the linear form on $i^G_P \delta$ given  by the absolutely converging integrals:
$$ v \mapsto \int_{U^-}  \langle \cc{v} (u^-), v(u^-) \rangle du^-, v \in i^G_P \delta$$
is an isomorphism. \eer   The restriction of functions to $K_0$  determines a bijection between $i^G_PE_\chi$ and $i^{K_0} _{K_0 \cap P} E$. If $v$ is an element of $i^{K_0} _{K_0 \cap P} E$, $v_\chi$ will denote its unique extension to an element of $i^G_PE_\chi$. 
\ber \label{vdel} Let  $\mm{V}(\delta, H)= \oplus_{x\in\mm{X}^G_M } \mm{V}(\delta,x,H) $  where   $\mm{V}(\delta,x, H) =(E')^{M\cap x^{-1} .H}$.\eer  Let $\eta= (\eta_x)_{x \in \mm{X}^G_M} \in \mm{V}(\delta,H)$.
  Let $J_\chi $ be the  subspace of  elements of $i^G_P E_\chi$  whose support  is contained in  $P \mm{W}^G_M H$ which is the union of the open $(P,H)$ double cosets in $G$.   One defines a linear form on $J_\chi$ by 
 $$ \langle \tilde{\xi}(P, \delta_\chi, \eta),  v\rangle = \sum_{x \in {\mm{X}^G_M} } \int _{M\cap x^{-1}.H\bb x^{-1}.H} \langle \eta_x , v ( yx^{-1})   \rangle dy  , v \in J_\chi.$$
From [BD], Theorem 2.7, one sees that 
  \ber \label{uniqueext}There exists   a non zero product $q$  of  functions on $X(M)_\si $ of the form  $ \chi\mapsto \chi(m) -c$,  for some $m\in M$ and $c\in \C^*$,  such that if  $q(\chi)\not=0$, $\tilde{\xi }({P, \delta_\chi, \eta} )$ extends to a unique $H$-invariant linear form on $i^G_P E_\chi$, denoted by ${\xi}({P, \delta_\chi, \eta}) $.
  Moreover for every $v$  element of $i^{K_0} _{K_0 \cap P} E$,   the map $\chi \mapsto q(\chi)\langle  \xi({P, \delta_\chi, \eta}, v_\chi \rangle$ extends to a polynomial function on $X(M)_\si$.  \eer
   
  When $ \xi({P, \delta_\chi}, \eta)$ is defined,  one defines   for $v \in i^G_P E_\chi$:
$$E(P, \delta_\chi, \eta,v) ( \dot{g} ) =\langle  \xi({P, \delta_\chi}, \eta),  (i^G_P\delta_\chi)  (g)v \rangle, g \in G. $$
Now, one uses ( \ref{bdl})  which  extends results of [BD] and [L]  when  the characteristic of $\F$  is  equal to zero  to the case where this characteristic  is  different from 2.  
From (\ref{bdl}),  [L],  Theorem 4 (ii),  [BD], Theorem 2.14 and Equation (2.33),  one sees that   if  $\chi \in X(M)_\si$ is such that $ Re( \chi \delta_P^{-1/2}) $ is strictly  $P$-dominant, $ \xi({P, \delta_\chi}, \eta)$ is defined and one has:
\beq \label{achtung} E(P, \delta_\chi, \eta, v) ( \dot{g} ) = \sum_{x \in \mm{X}^G_M} \int _{M\cap x^{-1}.H\bb x^{-1}.H} \langle \eta_x, v ( yx^{-1}g)  \rangle dy  , g \in G ,  Êv \in i^G_P E_\chi \eeq 
the integrals being convergent.


\subsection{Some results of Nathalie Lagier}
  One has the following  assertion which follows from [W], Theorem IV.1.1. 
Let $P=MU, P'=MU'$ be two $\si$-parabolic subgroups of $G$ with Levi subgroup $M$. 
\ber  \label{R} There exists $R>0$ such that  if  $\chi\in X(M)_\si$ satisfies  
$$  \langle {\rm Re} (\chi) , \aa \rangle >R, \aa \in \Delta(P) \cap \Delta(P'^-) , $$  the following integrals are convergent:
$$(A(P', P,\delta_\chi) v)(g):= \int_{U\cap U' \bb U'} v(u'g) du', v \in i^G_PE_\chi $$
Then $A(P', P, \delta_\chi)$ is an intertertwining operator between $i^G_P \delta_\chi $ and $i^G_{P'} \delta_\chi $. \eer 
The following results   are due to Nathalie Lagier (cf. [L], Theorem 5). We use the notation and hypothesis of the preceding subsection. 

Let $P$  be a standard $\si$-parabolic subgroup of $G$. Let  $(a_n) $  be a sequence in $A_M$ such that $(a_n) \to_{P}\infty$ i.e. such that for every root  $\aa$ of $A_M$ in the Lie algebra of $U$,  $(\vert \aa (a_n)\vert_F)$ tends to infinity.

 Let $(\delta,E)$ be a smooth unitary irreducible representation of $M$ and let  $\mu_{\delta}$ be its central character.    Let $\chi \in X(M)_\si$. Let us assume that  the real part of $\tilde{\chi}:=\chi\delta^{-1/2}_P $  is strictly $P$-dominant and satisfies (\ref{R})  for $P'=P^-$. Let $v \in i^G_P E_\chi  $  and $g \in G$. Recall that we have choosen $\mm{X}^G_MÊ\subset \mm{X}^G_{M_\vid}$ such that $\dot{1}\in \mm{X}^G_M$.  Then one has: 
\ber   \label{th5}If $\eta\in \mm{V}(\delta,x,H)$ with $x\in \mm{X}^G_M$   different from $\dot{1}$, one has:
$$lim_{n\to \infty}  \tilde{\chi}(a_n^{-1})\mu_{\delta}(a_n^{-1})   E(P, \delta_\chi, \eta, v) (\dot{1} a_ng)= 0.  $$\eer
and 
 \ber \label{th5ii}  If   $\eta\in \mm{V}(\delta,1, H)$, i.e. $\eta \in E'^{M\cap H}$, one has  the equality of 
$$lim_{n\to \infty} \tilde{\chi}(a_n^{-1})\mu_\delta(a_n^{-1}) E(P, \delta_\chi, \eta, v) (\dot{1} a_ng) $$
with $$\langle \eta, (A(P^-, P,\delta_\chi) v )(g)\rangle, $$   
\eer

Let $\ep $ be the trivial representation of $M_\vid$. Let $\chi \in X(M_\vid)_\si $ such that  the real part of $\tilde{\chi}:= \chi \delta_{P}^{ -1/2}  $ is  strictly $P_\vid$-dominant. Let $\eta$ be the linear form on $\C$ corresponding to $1$   and let $x\in \mm{X}^G_M$. We consider the  Eisenstein integrals for  $x^{-1}.H\bb G$. Then $ x^{-1}$ might be viewed has an element of a set $\mm{X}^G_M $ for $ x^{-1}. H$.  We view $\eta$ has an element of $E'^{M\cap H}= \mm{V}(\ep,1, H)$   and of $ E'^{ M\cap xx^{-1}. H}=\mm{V}(\ep, x^{-1}, x^{-1}.H)$.  Let $v\in i^G_{P_\vid} \chi$. We denote by $E_x (P_\vid, \eta, v)$ the Eisenstein integral for $x^{-1}. H \bb G$. 
Then one has:
$$E(P_\vid, \chi, \eta, v)(xg)= E_x (P_\vid,Ê\chi,  \eta, v)((x^{-1}. H) g), g \in G,$$ 
as it follows easily from  (\ref{achtung}). 
Using this, it follows from   [L],  Theorem 6: 
\ber  \label{th6vid}   Let $\ep $ be the trivial representation of $M_\vid$. Let $\chi \in X(M_\vid)_\si $ such that  the real part of $\tilde{\chi}:= \chi \delta_{P}^{ -1/2}  $ is  strictly $P_\vid$-dominant.  Let $P=MU$ be a   standard $\si$-parabolic subgroup of $G$. Let $(a_n)$ be a sequence in $A_M$ such that 
$(a_n ) \to_P \infty $. \\
Let $\eta$ be the linear form on $\C$ corresponding to $1$.  Let $x\in \mm{X}^G_M$. We view $\eta$ as an element  of  $\mm{V}(\ep,  1,H)$. For $v \in V:=i^G_P \C_\chi $, let $E_v :=E(P_\vid, \chi,\eta_x, v)$. Then  the sequence $(\tilde{\chi}(a_n^{-1}) E_v(x a_n) )$ has a limit. If $x\notin \dot{1} P$ this limit is equal to zero. Moerover if $x=\dot{1}$, one has: 
 $$\lim_{n \to \infty} (\tilde{\chi}(a_n^{-1}) E_v(\dot{1} a_n)= l(v), v \in V $$
 where $l$ is a non zero linear form on $V$.
 \eer
 Actually $l$ is explicit but what is important for us here is that it is non zero.
\subsection{ Applications of the results of N. Lagier }
 \begin{lem} \label{lem2} Let  $P=MU$ be  a standard  $\si$-parabolic subgroup of $G$. Let $(a_n)$ be a sequence in $A_M$ such that $( a_n) \to_P \infty$. If $(g_n) $ is a sequence in $G$ converging to $g \in G$ and such that for all $n\in \N$, $\dot{1} a_n g_n = \dot{1} a_n $,   then $g $ is an element of $H_P= U^-(M\cap H)$. 
\end{lem}
\dem  
One applies (\ref{th5ii}). We use the notation of this result.  If $J$ is a compact open subroup of $G$,  for $n$ large enough $g_n J= g J$. Hence, if $v \in i^G_P E_\chi$,
$$E(P, \delta_\chi, \eta, v)(\dot{1}a_n g_n ) = E( P, \delta_\chi, \eta, v)(\dot{1}a_ng),$$
for $n$ large enough.

First, let  $\delta$  be the trivial representation of $M$. One applies (\ref{th5ii}) to $v$ and $(i^G_P \chi(g)) v$  in  order to deduce from the preceding equality
$$(A(P^-, P, \chi) v) (g)=  (A(P^-, P, \chi) v)(1), v \in i^G_P \C_\chi$$
for $\chi $ sufficiently $P$-dominant. If $\chi $ is such that $A(P^-, P, \chi ) $ is bijective, one deduces the following equality:
\beq \label{eomega} v(g) =v(1), v \in i^G_{P^-} \C_\chi.\eeq 
Let us show that this implies   $g \in U^-M$. Let us write  $g  = p^-k $ with $k\in K_0$ and $p^-\in P^-$ . If $k \notin K_0\cap P^-$, there exists $v \in i^G_P \C_\chi $ such that $v(k)=0$ and $v(1)=1$, as the space of restrictions to $K_0$ of the elements of $i^G_{P^-} \C_\chi$ is equal to $i^{K_0}_{K_0\cap P^-} \C$.   This is  a contradiction to (\ref{eomega}).  Hence $g= u^-m $ with $u^-\in U^-$ and $m\in M$. \\Then applying (\ref{th5ii})   to any $(\delta, E, \eta)$,  we get similarly  the equality:   $$<\delta'(m)\eta, e>= <\eta, e>, e \in E.$$ The abstract Plancherel formula (cf. [Ber], section 0.2)  for $H\cap M \bb M$ implies $m \in  M\cap H$.  \qed

\begin{lem} \label{7} 
Let $P=MU, P'=M'U'$ be  two  standard  $\si$-parabolic subgroups of $G$. 
 Let $(a_n) $ (resp.,  $( a'_n)$) be a sequence in $A_M$ (resp. $A_{M'}$) such that $(a_n)\to_P \infty $ ( resp. $(a'_n)\to_{P'} \infty $). Let $J$ be a compact open subgroup of $G$. Let us assume that there exists $g, g' \in G$ such that for all $n\in \N$, 
$\dot{1} a_n g J = \dot{1} a'_n g' J$. Then, taking possibly subsequences, one has:\\
(i)   for all $\chi$  such that the real part of $\tilde{\chi}:=\chiÊ\delta_P^{-1/2}$ is  strictly $P_\vid$-dominant  $\tilde{\chi}Ê(a_n ^{-1} a'_n) $ has a non zero limit. \\
(ii) The sequence $(a_n ^{-1} a'_n) $ is bounded. \\
(iii) One has $P=P'$. \\
 (iv) If $Q$ is a $\si$-parabolic subgroup of $G$ such that $P\subset Q$, one has $\dot{1}_{Q} a_ng J =  \dot{1}_{Q}  a'_n g' J$ for $n$ large. 
\end{lem}
\dem (i) For all $n\in \N$, there exists $j_n\in J$ such that \beq \label{jn=}\dot{1} a_ng = \dot{1} a'_n g' j_n .\eeq  As $J$ is compact, one may take a subsequence and we may assume that $(j_n)$ converges to $j\in J$. Let $g''= g' jg^{-1}$. One will  apply  the result (\ref{th6vid}).   With its notations, let $v\in i^G_{P_\vid} \C_\chi $ and let us denote by $E_v$ the function $E(P_\vid , \chi, \eta,v)$. As $E_v $ is right invariant by an open compact subgroup of $G$,  one has $E_v ( \dot{1} a'_n g'j_ng^{-1} ) =  E_v( \dot{1} a_n'g'') $  for $n$ large. From (\ref{th6vid}), one has: \beq \label{limEv} lim_n \tilde{\chi} (a_n^{-1} ) E_v (\dot{1} a_n)= l(v), lim_n \tilde{\chi} ({a '_n}
^{-1} ) E_v (\dot{1} a'_ng'j_ng^{-1} )= l'(v) \eeq 
where $l$, $l'$  are  non zero linear forms on $i^G_{P_\vid} E_\chi$.
Also from (\ref{jn=}) one has:
\beq \label{Han} \dot{1} a_n = \dot{1} a'_ng'j_ng^{-1} \eeq
Let us show that there  exists $v_1\in V=i^G_{P_\vid} \C_\chi $ such that $l(v_1)$ and $l'( v_1)$ are non zero.  Let $v \in V$ such that $l(v)\not=0$. Then  $l $ does not vanish on  $v+{\rm Ker}(l)$. If  $l'$  vanished identically on  $v+{\rm Ker}(l )$  it would vanish on $V$, a contradiction which shows   that  $l' $ 
does not vanish identically on $v+{\rm Ker}(l )$. This proves our claim.

For such a $v_1$, one sees from (\ref{limEv}) and (\ref{Han})   that: \ber  \label{seqchi}The sequence $(\tilde{\chi} ( a_n{a'_n }^{-1}))$  tends to a non zero limit. \eer This proves (i).\\
(ii) By varying $\chi$ such that $\chi=Re \chi$  and such that $Re\chi $ describes a basis of $\a_\vid^*$  one gets  (ii). \\ (iii) 
If $P$ is different from $P'$,  by exchanging possibly the role of $P$ and $P'$,  there exists  a simple root $\aa$ of $A_\vid $   in the Lie algebra of  $U$ which is not a root in  the Lie algebra of  $U'$, hence which is a root in the Lie algebra of $M'$. Then $\vert \aa(a'_n)\vert_\F = 1$  and $ \vert  \aa (a_n{a'_n }^{-1})\vert_\F  $ is unbounded.  This would contradict (ii). Hence $P=P'$ and (iii) is proved.\\ 
 (iv)  From (ii), one writes 
$a'_n = a_n b_n$  where the sequence $(b_n)$ in $A_M$ is bounded.  Taking a subsequence we can assume that $(b_n)$ converges to $b\in A_M$.\\
Taking into account (\ref{jn=}), 
 one has $\dot{1}  a_n  = \dot{1} a_n c_n$ where $ c_n =  b_n g'j_n g^{-1}$. One deduces  from Lemma \ref{lem2} that the limit $c $ of $(c_n)$ is in $H_P$.  
 One has $a'_ng'J= a_nc_n gJ$. Hence for $n$ large one has:
$$\dot{1}_Pa'_ng'  J= \dot{1}_P a_n c_n g  J= \dot{1}_P a_n cg J $$
As $c  \in H_P$ and as $a_n\in A_M$ normalize $H_P$, one deduces that, for $n$ large:
$$\dot{1}_Pa'_n g'  J = \dot{1}_P a_n g J.$$
This proves (iv) for $Q=P$.

  Let $g, g'\in G$. In view of Theorem  \ref{theoct},  applied to the right transtlates  of $f$ by $g$, $g'$, there exists  $N\in \N$  such that for all $n\in \N$ greater than $N$ and for all $f\in C^{\infty}(H_Q\bb G)$ which is $J$-invariant $(c_{P,Q} f) (\dot{1}_P a_ng)= f(\dot{1}_Q a_n g)$ and $(c_{P,Q} f)   (\dot{1}_Pa'_ng')= f(\dot{1}_Q a'_ng')$.
Let $f$ be the characteristic function of $\dot{1}_Q a_n gJ\subset X_Q $. Let $n$ be an integer greater than $N$ and  let $x= a_ng$, $x'=a'_n g'$. From the above remark, one  has:
$$ (c_{P,Q} f) (\dot{1}_Px)=f(\dot{1}_Qx)=1,$$
$$(c_{P,Q} f)  (\dot{1}_P x')=f (\dot{1}_Q x').$$
  From (iv) for $Q=P$, one has  $\dot{1}_Px'J= \dot{1}_PxJ$. By $J$-invariance, this implies:  
$$(c_{P,Q} f) (\dot{1}_Px')= (c_{P,Q} f) (\dot{1}_P x).$$
Hence, one has 
$$f(\dot{1}_Qx')= 1 $$ and $x'\in \dot{1}_Q xJ$. This implies  $\dot{1}_QxJ= \dot{1}_Q x'J$. \qed 
\begin{lem} \label{8}Let $P=MU,P'=M'U'$ be two standard   $\si$-parabolic subgroups of $G$. Let  $(a_n)$ (resp., $(a'_n)$)
 be a sequence in $A_M$ (resp., $A_{M'}$) such that $(a_n)\to_P\infty $ (resp.,  $(a'_n)\to_{P'}\infty $). Let $g,g'\in G$ and $x, y \in \mm{X}^G_{M_\vid}$. Let us assume that the sequences 
$(x a_ng J)$ and $( y a'_ng' J)$ are equal. Then one has   $P=P' $,   $xP=yP$  and  $y=xm$  for some $m\in M$.
\end{lem}
\dem
Let $\chi \in X(M_\vid)_\si $ such that $\chi= \vert \chi\vert$ and such  that ${\rm Re}(\tilde{\chi})$ is strictly $P_\vid$-dominant. 
By exchanging possibly the role of  $x$ and $y$, and by taking a subsequence, one may assume that $\tilde{\chi}(a_n) \geq \tilde{\chi} (a'_n)$.  
Changing $H$ into $x^{-1}.H$,  one is reduced to the case where  $x=\dot{1}$.
Using the notation and the result of   (\ref{th6vid}), one sees that  there exists a non zero linear form $l$ on $V_\chi:= i^G_{P_\vid} \C_\chi$ such that for all  $v \in V_\chi$, one has: 
 \beq \label{lv} lim_n \tilde{\chi} (a_n)^{-1} E_v (\dot{1} a_ng ) = l (v).    \eeq

Let  $v\in V_\chi$ such that $l(v)\not= 0$.  One chooses $j_n\in J$ such that $\dot{1} a_ng= y a'_n g'j_n$. By extracting a subsequence, one may assume that $(j_n) $ converges to $j \in J$.  One has: \ber \label{Evl}$ E_v (\dot{1} a_ng)  = E_v ( ya'_ng'j)$ for $n$ large\eer  and $\tilde{\chi}(a_n) \geq \tilde{\chi} (a'_n)$.  
Let us assume  $y\notin \dot{1}P'$. Then from (\ref{th6vid})  
$$lim_n \tilde{\chi} (a'_n)^{-1}  E_v ( y a'_ng'j) = 0.$$ 
Together with (\ref{Evl}) this contradicts (\ref{lv}). Hence $y\in \dot{1} P'$ which implies (cf. (\ref{CD71}))  $y=\dot{1} m'$ for some $ m'\in M'$.   This implies the equality $y a'_n g'= \dot{1}a'_n m'g'$ as $a'_n \in A_{M'}$.  Hence one has $\dot{1} a_ngJ= \dot{1} a'_n m'g' J$. 
 Using  Lemma \ref{7}, one sees that $P=P'$. Hence $M=M'$ and the lemma follows. \qed 
  \section{ End of proof of Theorem \ref{theoexp}} \label{endtheo} 
  \setcounter{equation}{0}
\subsection{Definition of $exp_{X_P,X, J}$ \label{defexp}}
We have to deal only with the case $Q=G$ i.e. $X_Q=X$. 
If it does not exist a constant $C>0$ satisfying (i) of Theorem \ref{theoexp} for $Q=G$, there would exist  $x, y\in \mm{X}^G_{M_\vid}$,  and  sequences $(a_n), (a'_n) \in A_\vid$, $(\omega_n), (\omega'_n) \in \Omega$  such that $xP \not= yP$, $(\vert \aa ( a_n) \vert_\F$) tends to infinity for all roots $\aa$  of $A_\vid$ in the Lie algebra of $U$ and such that:
$$xa_n \omega_n J=  y a'_n \omega'_nJ, n\in \N.$$
By extracting subsequences, one may assume that $ \omega_n$ (resp., $\omega'_n$) converges to $ \omega$ (resp. $\omega'$). 
Let $Q=LV$ be the standard $\si$-parabolic subgroup of $G$ such that for $\aa\in \Delta(P_{\vid})$, the sequence $(\vert \aa(a_n)\vert _\F)$ is unbouded if and only if $\aa \in \Delta(Q,A_\vid)$. Clearly one has $Q\subset P$. 

By extracting subsequences, one will show that one can write     $a_n= b_n c_n$ where   the  sequence $(b_n) $ in $ A_{L}$ satisfies  $(b_n)Ê\to_{Q}\infty$ and where the sequence $(c_n)$ converges in $G$.  Let $( \delta_1, \dots, \delta_p)$ be  the union of $\Delta(Q, A_\vid)$ viewed as subset of $\a_\vid'$  and of a basis of $\a'_G$  viewed as a subset of $\a_\vid'$ (cf. (\ref{oplus})).  Let us look to the map $\phi: A_{\vid} \to  \R^p$  given by 
$a \mapsto (\delta_1(H_\vid(a)), \dots , \delta_p(H_\vid(a))$. Its image is a lattice of dimension $p$ as the image $\a_{\vid, \F}$ of $A_\vid$ by  $H_\vid$  is a  lattice of dimension equal to the dimension of $\a_\vid$. 
Its restriction to $A_L$ has the same property as it  factors through $H_{L}$ and $(\delta_1, \dots, \delta_p)$ might be viewed as a basis of $\a'_L$. 
Hence $ \phi(A_L)$ is of finite index in $ \phi (A_\vid)$. Hence one can find  $x_1, \dots,  x_q \in A_\vid $ such that for all $a \in A_\vid $ there exists $b\in A_L$ and 
$i\in \{1, \dots, q\}$ such that $\phi (a)= \phi(b x_i)$. This allows to define $b_n$ and $c_n = a_n(b_n)^{-1}$.  One has $c_n = x_{i_n} $ for some  $i_n\in \{1, \dots, q\}$. Then extracting a subsequence one may  even assume that $(c_n)$ is constant hence it converges. Moreover  as $\phi (b_n)= \phi (a_n)- \phi(x_{i_n})$ one has $ (b_n) \to_Q \infty $

Hence, for $n$ large,   $xa_n \omega_nJ =  x b_n c\omega  J$ where $c $ is the limit of $(c_n)$ . We introduce similarly $Q'$, $b'_n $ and $c'_n$. 
From Lemma  \ref{8} applied to $ G$ one deduces 
 $Q'=Q$  and  $HxQ=HyQ$. Hence, as  $Q\subset P$, one has  $xP= yP$. A contradiction which shows that there exists $C>0$ which satisfies (i). It is clear that any constant greater than such a constant enjoys the same property.\\
Let us assume that there is no constant satisfying (i) which satisfies also (ii). 
Proceeding as above, there would exist sequences $(b_n)$ in $A_L$ , $(b'_n)$ in $A_{LÔ}$, $c, c'\in G$, two standard $\si$-parabolic subgroups $Q=LV,Q'=L'V'\subset P$ of $G$ and  $x, y \in \mm{W}^G_{M_\vid} $ such that,  $(b_n)\to_{Q} \infty$, $(b'_n)\to_{Q'} \infty$,  and 
\ber \label{hyp1} $$xb_n c J = y b'_n c J.$$
 $$x_P b_n c J \not= y_P b'_n c J.$$ \eer
From Lemma \ref{8}, one sees that $Q=Q'$ and $ x\approx_Ly$. In particular $y= xl$ for some $l\in L$
 and, as $l$ commutes to the elements $b'_n$ of $A_L$,   one has:
 $$xb_n c J = xb'_nl c' J.$$
Conjugating  $x^{-1}$, one gets an  equality of left $x^{-1}. H$ cosets.  From Lemma \ref{7} (i), applied to $x^{-1}. H$ instead of $H$, one deduces that $(b'_n{ b_n} ^{-1})$ is bounded. Hence, by  taking a subsequence  one can assume that it has a limit.
Then from Lemma \ref{7} (iii)  one gets for $n$ large:
$$(x^{-1}.H)_P b_n c J= (x^{-1}.H)_P b'_n lc' J .$$
Hence there exists a sequence in $J$, $(j_n)$ such that 
$$(x^{-1}.H)_P b_n c j_n= (x^{-1}.H)_P b'_n l c'. $$
Hence  $b_n  c j_n  c'^{-1}l^{-1} b'_n  \in (x^{-1}.H)_P$.
As the stabilizer of  $x_P$ is equal to $(x^{-1}.H)_P$, one deduces from this  the equality:
$$x_Pb_nc j_n=x_P b'_nlc  $$
As  $x = y l$ and $l\in L \subset M$, one has $x_P=y_P l$. As $ l \in L$ commutes to $b'_n  \in A_L$, 
one deduces from this the equality $$x_Pb_nc J=y_P b'_nc'J,  $$ for $n$ large.
This  contradicts our hypothesis (\ref{hyp1}). Hence there exists $C>0$ which satisfies (i) and (ii).

\subsection{Injectivity of $exp_{X_P, X,J}$}
 Let us prove that one can choose $C>0$ such that $exp_{X_P,X}$ is injective on $N_{X_Q , J } (P,C)$.   
Let us assume that  every constant  $C>0$ satisfying conditions (i), (ii)  of Theorem \ref{theoexp} does not satisfy condition (iii). From the finiteness of $\mm{X}^G_{M_\vid}$   and proceeding  as in  section \ref{defexp}, one sees that  there would exist   $ x, x' \in \mm{X}^G_{M_\vid} $,   two $\si$-parabolic subgroups $Q=LV, Q'=L'V'\subset P$ of $G$,   a sequence $ (a_n) $ in $A_L$, a sequence $(a'_n)$ in $A_{L'} $ such that $(a_n)\to _{Q} \infty$,  $(a'_n)\to _{Q'} \infty$ and two elements $d$ and $d'$  of $A_\vid\Omega$ such that:  $$x a_n dJ \not=  x' a'_n d'J$$
and $$x_P  a_n d J =  x'_P  a'_n d' J.$$ Let $ f_n$ be the characteristic function   of $x a_n d J $.   For  $n_0  $ large enough one  can use  Theorem  \ref{theoct} (iv) for the right translates of $f_{n_0}$ by $d$ and $d'$ and one has, by setting $a=a _{n_0}, f=f_{n_0},$ etc.:

$$f (x a d )  =(c _{P, G} f )(x_P a d ), f( x' a' d' )= (c _{P, G} f ) (x_P' a' d' )$$
But, by our assumptions $f (x a d)= 1$ and $x_P  a d J =  x'_P a' d' J$. Hence, by $J$ invariance, one has:
$$f( x 'a' d' ) = 1$$
which implies $$ x  a d J =  x'  a' d' J.$$
This is a contradiction to our hypothesis.  This achieves to prove that there exists a constant $C>0$ such that the properties (i), (ii) and (iii) of Theorem \ref{theoexp} are  satisfied. 
\subsection{Volumes}
 The following lemma will allow to finish the proof of Theorem \ref{theoexp}.
\begin{lem} \label{lem1}Let $K$ be a compact open subgroup of $G$ with a strong  $\si$-factorization for  standard $\si$-parabolic subgroups (cf. (\ref{sifac})). Let $P=MU$ be a standard $\si$-parabolic subgroup of $G$. Let  $a \in A_\vid$ which is $P_\vid$-dominant. Then \\ (i)$$Ha K=  H a K_M K_U.$$
where $K_M=K\cap M, K_U= K \cap U$. \\
(ii)  $$vol_X  \dot{1}  a K=  vol_{X_P} \dot{1}_P a K.$$ \end{lem}
\dem
(i) As $ K_{M_\vid} K_{U_\vid} = K \cap P_\vid$  and $K_M K_U =K \cap P$, it is enough to prove  (i) when $P= P_\vid$. Let us assume this in the sequel. 
If $u^- \in K_{U^-} $, as $a$ is $P_\vid$-dominant,  one has $a.u^-= a u^{-} a^{-1} \in K_{U^-} \subset K= K_H K_M K_U$ (cf. (\ref{sifac}) (ii) and  (iii)).  
Hence one has: $$ H au^-= H(a.u^- ) a  \in H K_M K_U a.$$
But $K_MK_U a= a(a^{-1}. K_M) (a^{-1}. K_U)$. As $M=M_\vid$ and $a \in A_\vid$, $a^{-1}. K_M= K_M$. As $a $ is $P_\vid $-dominant $a^{-1}. K_U \subset K_U$ (cf. (\ref{sifac}) (ii)).  Altogether, this shows: $$ H a K_{U^-} \subset Ha K_ {M} K_{U}$$ 
One deduces (i) from the equality $K=K_{U^-}K_MK_U$. 
\\ Let us prove (ii).  Let $P$ be a standard $\si$-parabolic subgroup of $G$.  As $U^- \subset H_P$ and $K=K_{U^-}K_MK_U$,   and $a.K_{U^-} \subset {U^-}$,  (cf. (\ref{sifac}) (ii) ) one has: $$ \dot{1}_P a K=\dot{1}_P a K_MK_U . $$ Then (ii) follows from (i), from  the fact that $a K_MK_U \subset  P$ and from our choice of measure on $X_P$ (cf. section \ref{mea}).  \qed \\
{\em End of proof of Theorem \ref{theoexp}}
\\ Let $K$ and $J$ as in the theorem. The proof of (iv)  reduces to prove the statement for subsets of $N(x,P, C)/J$ for $x \in \mm{X}^G_M$.  Using our choices of volumes and translating  sets on the left  by $x^{-1}$ and changing $H$ in $x^{-1}.H$, one is reduced to the case $x=1$.  For  
 $K$  and $J$ as in the theorem, we have:
  $$\omega.J \subset K, \omega\in \Omega.$$
  Let $\omega\in \Omega$ and one sets    $J':=\omega.J\subset K$. As $\Omega$ is compact and is left $K$-invariant,  $\Omega/J $ is finite and   $J'$ varies in a finite set.  Let us assume that  $C>0$ satisfies Theorem \ref{theoexp} (i), (ii) and (iii) for  all groups $J'$ .  One has to prove that for $a\in A_\vid^+(P,C)$:
 $$  vol_X (\dot{1} a\omega J)= vol_{X_P} (\dot{1}_Pa \omega J).$$
 As the measures on $X$ and $X_P$ are right invariant by $G$, in order to prove this  equality, 
 it is enough to prove the equality:
 $$vol_X (\dot{1} aJ')= vol_{X_P} (\dot{1}_Pa J').$$
 Let $K_a$ (resp.,   $K'_a$) be the stabilizer in $K$ of $\dot{1}a$ (resp., $ \dot{1}_P a$). 
 We need the following fact. Let $K_1$ be a closed subgroup of $K$. Let us assume that a Haar measure is given on $K$ and let $K_1\bb K $ be  endowed with  the image of this measure. Let $ X\subset K$ and $Y$ its  image in $K_1 \bb K$ . Then  $vol_{K_1\bb K} (Y)= vol_K (K_1X)$. 
 From this applied to $K_1= K_a$ and $K_1= K'_a$ and from Lemma \ref{lem1} (ii), it is enough to prove the equality:
 $$K_a J'= K'_aJ'.$$ 
  The image of  the set  $\dot{1} a K'_a J'$   by the map  $exp_{X_P,X,J'}$   is equal   $\dot{1}_P a J'$, as it follows from the definition in Theorem \ref{theoexp} and the equality $\dot{1}_P a K'_a J'= \dot{1}_P aJ'$. From the definition of  $exp_{X_P,X,J'}$, this image is also equal to the image of $ \dot{1} aJ'$. Hence from the  part (iii) of Theorem  \ref{theoexp}, one deduces the equality: $$\dot{1}  aJ' =\dot{1}a K'_a J'$$
  Looking to the orbit of $\dot{1}a$ under $K$ one deduces from this the inclusion: 
  $$ K'_a J'\subset K_aJ' .$$
  We recall that $K \subset \Omega$. 
To prove the reverse inclusion let us remark  that  $\dot{1} a K_a  J' $  is   equal to   $\dot{1} aJ'$.  From the definition of $exp_{X_P, X, J'}$ 
one deduces the equality:
$$\dot{1 }_P a  K_a J'= \dot{1}_P a J' $$
which implies as above: 
$$ K_aJ'\subset  K'_a J'. $$
This implies  the required equality. This finishes the proof of of the theorem.\qed
  \section{Bernstein maps and Scattering Theorem}
  \setcounter{equation}{0}
 \subsection{Constant term and $exp$-mappings}
 The following proposition is an immediate corollary of Theorems \ref{theoct} and \ref{theoexp}.
 \begin{prop} \label{ctexp}
  Let $P\subset Q$  be two standard $\si$-parabolic subgroups of $G$. Let $J$ be a compact open subgroup of $G$ small enough to satisfy the conditions of Theorem \ref{theoexp}.
 There exists $C>0$  such that $exp_{X_P,X_Q,J}$ is well defined on 
 $N_{X_Q, J}(P,C)$ and satisfies for all $J$-invariant function $f$ on $X_Q$:
 $$(c_{P,Q}f) (exp_{P,Q,J}(xJ))= f(xJ), xJ \in N_{X_Q, J}(P,C).$$.
 \end{prop}
 
 \begin{rem}
 In [SV],  for $G$-split and $X$ spherical, the $exp$-mappings are  introduced before the maps $c_{P,Q}$,  by means of wonderful compactifications, and the maps $c_{P,Q}$ are  defined by the relation above.
 \end{rem}
 \subsection{Bernstein maps $e_{Q, P}$}
 We thank Joseph Bernstein for having suggested  to us the proof  of the following Theorem.
 \begin{theo} \label{theoe}Let $P=MU \subset Q=LV$ two standard $\si$-parabolic subgroups of $G$. The right $G$-invariant measure  on $X_P$ allows to identify  $C_c^{\infty}(X_P)$ to a subspace of the dual of $C^{\infty}(X_P)$. Let $e_{Q,P} $ be the  restriction of the transpose map of $c_{P,Q}$ to $C_c^{\infty}(X_P) $. \\ Let   $J$   and  let   $C>0$  be as in Theorem \ref{theoexp} . \\
 (i) Let $xJ \in N_{X_Q, J}(P,C)$ and $y =exp_{X_P,X_Q,J}(xJ)$. Then the image by $e_{Q,P}$  of the characteristic function  of $yJ\subset X_P$ is the characteristic function of $xJ\subset X_Q $. 
 \\(ii)  For $f\in C_c^\infty(X_P)$ supported in  $exp_{X_P,X_Q,J} (N_{X_Q, J}(P,C))$, $e_{Q,P} f $ has its support in $N_{X_Q, J}(P,C)$ and 
 $$ (e_{Q,P}f) (xJ))= f (exp_{X_P,X_Q} (xJ ) ),  xJ\in N_{X_Q, J}(P,C).$$ÊÊ
 (iii)  The map $e_{Q,P} $ has its image in $C_c^{\infty}(X_Q)^J$.   
 \end{theo}  \dem
  (i) 
 We fix a compact open subgroup $J$  and $C$ as in the preceding proposition from which we use the notations.
 Let $xJ \in N_{X_Q, J}(P,C)\subset X_Q/J$  .  Let $f$ be the characteristic function of $exp_{X_P,X_Q,J}(xJ) $ which is a  $J$-invariant function on $X_P$. Let $g \in C^{\infty} (X_Q)$. 
One has $$\langle e_{Q,P} f , g\rangle= \langle f, c_{P,Q}g\rangle$$
and by the preceding proposition one sees:
$$\langle e_{Q,P} f , g\rangle= g(xJ).$$
This implies that $e_{Q,P} f $ is the characteristic function of $xJ$. This proves (i).
\\  (ii)  follows by linear  combinations. 
\\ (iii)
Let $a\in A_M$ be strictly $P$-dominant. Let $y\in X_P$. From  the Cartan decomposition for $X_P$ one sees that for $n$ large, $a^ny $  is  of the form $a^ny=exp_{X_P,X_Q,J}(x_nJ)\in X_P/J$ for some $x_nJ \in N_{X_Q, J}(P,C)\subset X_Q/J$.

For $n\in \Z$,  let $f_n$ be the characteristic function of $ a^nyJ \subset X_P$. One has just seen that for $n$ large in $\N$,  $e_{Q,P}(f_n) $ is in $C_c^{\infty}(X)^J$. Let us assume that it is not true for all $n\in \N$.  Then there would exists $N\in \N $ such that $e_{Q,P}(f_n)\in  C_c^{\infty}(X)^J $  for $n > N$ and such that $e_{Q,P}(f_N)\notin  C_c^{\infty}(X)^J $.

We want to apply Theorem  A of [AAG] in order to prove  that the $C_c^\infty(G)^J$-module $C_c^{\infty}(X_P)^J $  is finitely generated.
For this it is necessary to see that one may apply it to each homogeneous space $x_PG$ which is isomorphic to $U^-(M\cap x^{-1}.H)Ê\bb G$. The first thing to prove is that for each parabolic subgroup $R$ of $G$, the number of $ (U^- (M\cap x^{-1}.H),R)$-double cosets is  finite. By using conjugacy, one can assume that $R$ contains $A_0$. By the Bruhat decomposition,  one has  $G= \cup_{i} Px_i R$, where  $(x_i)$ is a finite family  of elements of $G$ normalizing $A_0$. It is enough, to prove our claim, to show that for each $i$, $R_i:=( x_i.R )\cap M$ has a finite number of orbits in the symmetric space $(M\cap x^{-1}. H) \bb M$.  But  $R_i$ is a parabolic subgroup  of $L$ and our claim follows from  [HW], Corollary 6.16. 

The second thing to prove, in order to apply Theorem A of [AAG] is that : \ber \label{finiteinv} For each  finite length smooth $G$-module $V$, the dimension of the space $V'^{U^-(M\cap x^{-1}.H)}$ is finite. Ê\eer
But this dimension is precisely the dimension of $j(V)'^{M\cap x^{-1}. H)}$ where $j(V)$ is the Jacquet module of $V$ with respect to $P^-$. This space is finite dimensional (cf. [D], Theorem 4.4.)

Now, one can apply Theorem A of [AAG] to conclude that the $C_c^\infty(G)^J$-module $C_c^{\infty}(X_P)^J $  is finitely generated.
  Moreover the algebra $C_c^\infty(G)^J$
 is Noetherian (cf. [R] Corollary of Theorem VI.10.4). \\ Hence, it follows that an ascending chain of $C_c^\infty(G)^J$-submodules  of $C_c^{\infty}(X_P)^J $  is stationnary.\\
We apply this to the $C_c^\infty(G)^J$-submodules of $C_c^{\infty}(X_P)^J $, $M_n$, generated by $f_0, \dots f_{-n}$. 
Hence there exists $n\in \N$ and $ \phi_0, \dots \phi_n \in C_c^\infty(G)^J$ such that:
$$f_{-n-1}=f_0*\phi_0+Ê\dots + f_{-n}* \phi_ n$$ 
 
 Using that the right $G$-action and the left $A_M$-action commute (cf. Definition \ref{aaction})   and applying the left  action of  $a ^{n+1+N}$ to  the above identity, one gets:
 $$ f_N= f_{n+1+N}*\phi_0+ \dots + f_{1+N} *\phi_{n}$$
From Theorem \ref{theoct},  $c_{P,Q}$ is a morphism of $G$-modules. Hence it is also the case for $e_{Q,P}$.   Hence $e_{Q,P}  (f_N)$ is in $C_c^{\infty}(X_Q)^J $. From the definition of $N$, we get a contradiction.  Hence in particular, $e_{Q,P}f$ is in  $C_c^{\infty}(X_Q)^J $. The theorem follows by linearity. \qed\\
Let $(\pi, V)$ be a smooth representation  of a parabolic subgroup $P=MU$ of $G$. One
denotes by
  $(\pi_P,V_P)$ the tensor product of the   quotient of $V$ by the
$M$-submodule generated by the $\pi(u)v-v$, $u\in U, v\in V$,  with the representation of $M$ on $\C$ given by $\delta_P^{-1/2}$.  We call it the
normalized  Jacquet module of $V$ along $P$.  We denote the natural projection map from $V$ to $V_P$ by
$j_P$  and sometimes $ \pi_P$ will be denoted $ j_P(\pi)$.
\begin{lem} \label{fpegal1} 
Let $P$ be  a semistandard $\si$-parabolic subgroup of $G$.\\
  (i)If $f\in C_c^\infty(X)$ has its support in $X_M P$ we define, using (\ref{dptriv}),  $f^{P} \in C^\infty (X_M) $ by 
$$f^P(xm)=\delta_P^{1/2} (m) \int_U f(\dot{x} m  u) du, x \in \mm{X}^G_M, m\in M $$
Then $f^P\in C_c^\infty ( X_M)$. \\
(ii)  The map $f\mapsto f^P$  goes through  the quotient to an intertwining map between the normalized Jacquet module $C_c^\infty(X_MP)_P$ of the $P$-module  $C_c^\infty(X_MP)$ 
 and $C_c^\infty(X_M)$.\\
 (iii) This intertwining map is bijective and its inverse define an intertwining injective map $m^X_P: C_c^\infty(X_M) \to C_c^\infty(X)_P$.\\
 (iv) One can replace $X$ by $X_P$ in (i), (ii) and (iii)  and one gets an injective intertwining map $m_P: C_c^\infty(X_M) \to C_c^\infty(X_P)_P$
\end{lem}
\dem
(i) follows  easily from the definition.\\
(ii) It is clear that our map goes through the quotient to a map between the normalized Jacquet module $C_c^\infty(X_MP)_P$ of the $P$-module  $C_c^\infty(X_MP)$.  
On the other hand, for $f \in C_c^\infty(X_MP)$ one has:
$$Ê(\rho(m_0) f)^P(xm)=\delta_P^{1/2}(m)\int_U f(xmm_0m_0^{-1}u m_0) du $$
One makes the change of variable $u'=m_0^{-1}u m_0$ to achieve to prove  the intertwining property of (ii).
\\ As an $U$-space,  $X_M P$ is isomorphic to $X_M\times U$ where $U$ acts trivially on the first factor. This implies easily (iii).\\
(iv) is proved similarly. \qed

\begin{prop} 
  We denote  by $j_P(e_P)$ the map between the normalized Jacquet modules $C^\infty_c(X_P)_P$ and $C^\infty_c(X)_P$ determined by $e_P:= e_{G, P} $.
Then $$j_P(e_P)\circ m_P = m_P^X.$$
\end{prop} 
\dem
One has to prove;
\beq \label{emm} j_P (e_P) (m_P(f))= m_P^X (f)\eeq 
for all $f \in C_c^\infty(xM)$ and $x \in  \mm{X}^G_M$.
 Changing $H $ to $x^{-1}.H$, one is reduced to prove  (\ref{emm}) for $x=1$.
One writes the Cartan decomposition for $M\cap H\bb M$:
$$M\cap H \bb M= \cup_{x \in \mm{X}^M_{M_\vid}}xA^+_\vid (P_\vid, P, 0) \Omega_M,$$ where $\Omega_M$ is a compact set of $M$ and $\mm{X}^M_{M_\vid}$ is the analog of $\mm{X}^G_{M_\vid}$. 
The $M$-module of  compactly supported smooth functions on $\dot{1}M$  is the linear span  the characteristic functions of $\dot{1}x a \omega J$ where $J$ describes  a basis of neighborhood of 1 in $M$ made   of compact open subgroup of $M$, $x \in \mm{X}^M_{M_\vid}$, $\omega \in \Omega_M$, $a \in A^+_\vid (P_\vid, P, 0)$ . As $m_P, m^X_P, e_P$ are $M$-equivariant, one  has to prove (\ref{emm})  for every $f$ among  a set of generators of this $M$-module  
Again we reduce to $x=1$. Taking into account  (\ref{aina}), one can  write  $a= a'b$ with $a'  \in A^+_\vid$. 
and $b \in A_M$. As $b$ commutes to $J$, one is reduced to prove (\ref{emm}) for the characteristic functions of $\dot{1} a \omega J$,  with $a \in A_\vid^+$ and $\omega\in \Omega_M$.

As $\omega J= \omega J\omega^{-1}\omega$, the characteristic functions 
of $\dot{1} aJ'$ where $J'$ describes the set of $\omega. J$ for  $J$ as above, $a\in A_\vid^+$, $\omega \in \Omega_M$ is a set of generators of $C_c^\infty (\dot{1}  M)$. 

Let $(J'_n) $ be as in Lemma \ref{strongsi}. By continuity and compacity, there exists a neighborhood  $\mm{V}$of 1 in $M$ such that:
$$\omega. \mm{V} \subset (J'_0)_M , \omega \in \Omega_M$$
One can assume that all the groups $J$ above are contained in $\mm{V}$. Hence all the groups $J'$ are contained in $(J'_0)_M$. For such a group, let $n\in \N$ such that $(J'_n)_M \subset J'$. 
Then as $J'$ is the disjoint union  of the left $(J'_n)_M$-cosets, the characteristic function of 
$\dot{1} a J'$ is a linear combination of the  characteristic functions  of $\dot{1} a j'(J'_n)_M$ where $j'$ describes $J'$. But as $J'_n$ is normal in $J'_0$ ( cf. Lemma \ref{strongsi} ) and $J' \subset (J'_0)_M$,  $(J'_n)_M$ is normal in $J'$. Hence $\dot{1} a j'J'_n= \dot{1} a J'_nj'$.
Hence, again by $M$-equivariance, one has to prove (\ref{emm})  for $f$ equal to the characteristic function of  $\dot{1} a (J'_n)_M$,  $n \inÊ\N, a \in A_\vid^+$. 

For simplicity we write $J$ instead of $J'_n$ and let $g= vol (J_U) \delta_P(a)^{1/2} 1_{\dot{1} a J_M}$ 
 and let $f=  1_{\dot{1}_P a J_M J_U} \in C_c \infty(X_P)$. Then $f^P = g$.
 Then, by definition of $m_P$, one has:   $$m_P(g)= j_P(f)$$ where $j_P(f)$ is the image  of  $f$ in the normalized Jacquet module of $C_c^\infty(X_P)$.  Similarly the characteristic function $h $ of  $\dot{1} a J_M J_U$  satisfies $h^P= g$. Hence one has: $$m_P^X(g)= j_P(h)$$ and   
$$(j_P(e_P)) ( m_P (g))= j_P(e_P (f)).$$
It remains to prove:
$$(j_P (e_P)) ( m_P (g)) =  m_P^X(g)$$ 
i.e.  $$(j_P (e_P)) (j_P(f))= j_P (h)$$
For this, it is enough to prove:
$$e_P (f)=h.$$
One has  $$\dot{1} a J_M J_U= \dot{1} a J$$ from Lemma  \ref{lem1}.
As $ J_{U^-}$ is normalized by $a \in A_\vid^+$  (cf. Lemma \ref{strongsi}),  one has $$   \dot{1}_P a J_M J_U= \dot{1_P} a J$$
Then the required equality follows from Theorem \ref{theoe}  (i). \qed

\subsection{Discrete spectrum}
An irreducible  subrepresentation of $C^\infty (X)$, $(\pi,V)$,  is said discrete if the action of $A_G$ is unitary and the elements of $V$ are square integrable mod $A_G$.  Obviously if $ \psi$ is element of the group  $X(G)_{\si,u}$ of unitary elements of $X(G)_\si$, the representation $\pi_\psi$ of $G$ in the space $V_\psi:= \{\psi v\vert v \in V\}$ is a also a discrete series. Moreover $\pi_\psi$ is isomorphic to $\pi\otimes \psi$. 
Let $\chi$ be a unitary character of $A_G$ and 
let $L^2(X, \chi)_{disc}$ the sum of all $X$-discrete series on which $A_G$ acts by $\chi$.
\begin{theo} \label{theodisc}  Let $J$ be a compact open subgroup of $G$ and $ \chi$ a unitary character of $A_G$. Then the space $L^2(X, \chi)_{disc}^J$ of $J$-invariants of $L^2(X, \chi)_{disc}$ is finite dimensional.
\end{theo} 
\dem Ä
One will see that the proof of  Theorem 9.2.1 of [SV] adapts by changing $Z(G)^0$ to $A_G$, and, for  a standard $\si$-parabolic subgroup $P=MU$ of $G$ by changing    $Z(X_P)$ to $A_M$ acting on the left.  . 

Let   $A^+_P$ be  the set of  $P$-dominant elements of $A_P$. 
 Let $N'_P$   be equal to $N_{X,J}(P,C)$ for $C>0$ large enough in such a way that the $exp$-maps are defined and such that the identity of Proposition \ref{ctexp} holds. Let $N_P= N'_P \setminus_{Q\subset P, Q \in \mm{P}, Q\not=P}   N'_{Q}$. Then  the $N_P$ covers $X$. We remark that  $exp_{X_P,X,J}(N'_P)$  is stable by the left action of $A^+_P$ as well as  $N''_P:= exp_{X_P,X,J}(N_P)$. One sees from the definitions that there is a finite  subset $\Omega_P$ of $X_P/J $, such that $ N''_P=A^+_P \Omega_P$.   Let $(\hat{A})_\C ^{J_M}$ be the set of complex  characters of $A_M$ which are trivial on $A_M\cap J$.  Let $\mm{P}$ be the set of standard $\si$-parabolic subgroups of $G$. We choose a map $R:\mm{P}\to \N$, $P\mapsto r_P$ and we define $\mathfrak{S}_{R}:=  \prod_{P\in \mm{P}}( (\hat{A})_\C ^{J_M})^{r_P}$. An element of $x \in \mathfrak{S}_R$ is denoted $[(\chi_{i})_{i=1, \dots,r_P}]_{P \in \mm{P}
}$. We consider for $a \in A_M$,  \beq \label{con}\prod_{i=1,\dots, r_P}( \mm{L}_a-\chi_{i}(a ))\eeq 
Let $x \in\mathfrak{S}_R$. We consider the   subspace $V_x\subset C^{\infty}(X)^J$ of $J$-invariant functions on $X$,  $f$,  such that for all  standard $\si$-parabolic subgroup $P$ of $G$  and $a \in A_M$,   $c_{P, G}f$ is annihilated by (\ref{con}). Then $V_x$ is invariant by the Hecke algebra  of $C_c^\infty $ functions on $G$ which are right and left invariant by $G$:  this is due to the fact that $c_{P, G}$ is a $G$-morphism and that the right  action of $G$ on $C^\infty(X_P)$ commutes with the left  action of $A_P$.

 Recall that from our hypothesis on $C$ that:$$(c_{P,G}f) (exp_{X_P, X,J} (x) = f(x), x\in N_P. $$
 Then $V_x$ is finite dimensional, as it is shown in the proof of Theorem 9.2.1 of [SV].
 The rest of the proof is entirely analogous to the proof of this Theorem. \qed
\begin{cor} \label{cordisc} Let $J$ be a compact open subgroup of $G$. 
There exists finetely many discrete series for $X$, $(\pi_i, V_i)$, $i=1, \dots, n$ such that any discrete series, $(\pi,V)$  for $X$ is of the form $(\pi_i)_\chi$ where $\chi$ is element of the group  $ X(G)_{\si,u}  $ of unitary elements of $X(G)_\si$ and $i\in \{1, \dots,n\}$.
\end{cor}\dem
Looking to Lie algebras one sees that restriction map from  the group $X(G)_{\si,u}$  of unitary elements  of $ X(G)_\si$ to the group $X(A_G)_u$ of  unitary elements of $ X(A_G)$ is surjective. On the other hand the action by multiplication  of $X(A_G)_u$ on $(\hat{A}_G)_u ^J$ has finitely many orbits (cf. \ref{surjection}). Hence one is reduced to the case where the restriction of the central character of $\pi$ is one of the representatives of these orbits. Then the corollary follows immediately from the Theorem.\qed 
The proof of the following Lemma is immediate.

 \begin{lem} \label{indcc}Let $\delta_{P, \mm{X}^G_M}$ be the function on $X_M$  such that,  for all $x \in \mm{W}^G_M $,  its restriction to $xM$ is equal to the function  $\delta_{P,x}$ occuring in (\ref{dxm}).  For a function $f$ on $X_P$ we associate  the map $T(f)$ on $G$ with values in the space of functions on $X_M$ defined by:
 $$(T(f)(g))(x)= \delta_{P, \mm{W}^G_M}^{-1/2}(x) f(xg), x \in X_M, g \in G$$
 (i) One has
 $$T(f)(mg) =  (\rho \otimes \delta_P^{1/2}) (m) f(g), m \in M, g \in G.$$
 (ii) The map $T$ induces  a bijective $G$-intertwining map  between  $C_c^\infty (X_P)$ and $ i^G_{P^-} C_c^\infty (X_M)$ (resp., $C^\infty (X_P) $ and $i^G_{P^-}C^\infty(X_M)$).
\\(iii) Let $\chi $ be a unitary character of $A_M$.  The map $T$ induces  a bijective isometric  $G$-intertwining map  between   $L^2(X_P)$ and the unitarily induced representation from $P^-$ to $G$ of $L^2(X_M)$ ( resp.,  $L^2(X_P, \chi)_{disc}$ and the unitarily induced representation from $P^-$ to $G$ of $L^2(X_M, \chi)_{disc}$). 
  \end{lem}
  \dem (i) is immediate.
  \\(ii) From (i), it remains only to prove the bijectivity. The inverse map to $T$ is easily described using the fact that $X_P=X_M\times _{P^-}  G$.
  \\(iii) follows easily from the definition of the scalar product on unitary induced representations from $P$ to $G$ (cf. (\ref{duali})) and from the definition of the $M$-invariant measure on $X_M$ (cf. (\ref{dxm}) and (\ref{dmx})). \qed 
 \begin{lem}  \label{inddisc}
 $L^2 (X_P)_{disc}$ is unitarily equivalent to the unitary induced representation from $P^-$ to $G$ of $ (L^2 (X_M)_{disc} )$
\end{lem}
\dem
 The Lemma  follows from the analog of Corollary 9.3.4 in [SV] and of  of Lemma \ref{indcc} (iii). Notice that this Corollary follows from
l.c. Equation (9.1). To establish its analog, one remarks that $A_M$ acts freely on the left on $X_P$. \qed 
\begin{lem} \label{dsconj} 
The $G$-space $X_P$ satisfies the discrete series conjecture 9.4.6 of [SV] for the parabolic subgroup $P^-
 $ and the  torus of unitary unramified characters of $P^-$, $D^*:=   X(M)_{\si,u} $.
\end{lem}
\dem From Corollary \ref{cordisc} of Theorem \ref{theodisc}, there is a denumerable family of $X(M)_{\si,u}$-orbits of discrete series. Then the Lemma follows from Lemma  \ref{inddisc}. \qed \subsection{Bernstein maps}
The proof of the following theorem is entirely analogous to the proof of Theorem 11.1.2 in [SV].  
\begin{theo} \label{Bmap}For every pair of standard $\si$-parabolic subgroups of $G$, $P\subset Q$, there exists a canonical $G$-equivariant  map  $i_{P,Q}: L^2(X_P) \to L^2(X_Q)$ characterized by the property that for any $\Psi \in C_c^{\infty} (X_P)  $ and any $a$ element of the set $ A_P^{++}$ of strictly $P$-dominant elements of $A_P$,  we have:
$$lim_{n \to \infty} (i_{Q,P}\mm{L}_{a^n}\Psi  - e_{Q,P} \mm{L}_{a^n} \Psi)=0$$
where the limit is in $L^2(X_Q)$.
\end{theo}
 
Then  as a corollary of Theorem  \ref{Bmap}  and of the analog of Proposition 11.6.1 of [SV], one has  the following  analog of l.c Corollary  11.6.2. The proof requires the criteria for discrete series of symmetric spaces due to Kato and Takano [KT2]:
\begin{prop} Let $L^2(X)_P$ the image of $L^2(X_P)_{disc} $ under $i_P:=i_{G,P}$. Then one has:
$$L^2(X) = \sum_{P\in \mm{P}_{st} }L^2(X)_P.$$ 
\end{prop}
\subsection{Scattering theory }
From Lemma  \ref{dsconj}, one  proves the analogous  of Proposition 13.2.1 in [SV]  in which we use $A_M$ and $A_L$ instead of $A_{X, \Theta}$ and $A_{X, \Omega}$ and where $P=MU$, $Q=LV$ are $\si$-parabolic subgroups of $G$.  This is a step for the analogous of Proposition 13.3.1 in l.c. . We will only recall part  (2) of it. 
\begin{prop} 
Let $P=MU,Q=LV$ two standard $\si$-parabolic subgroups of $G$. If the dimensions of $A_M$ and $A_L$ are distinct, $L^2(X)_P$ is orthogonal to $L^2(X)_Q$.
\end{prop} Let $\Theta_P$ (resp.,  $\Theta_Q$) be the set of elements of $\Sigma(P_\vid)$  which are trivial on $A_M$ (resp. $A_L$). 
    We define $W(P,Q)$ as the set of elements of $w\in W(A_\vid)$ such that $w(\Theta_P)= \Theta_Q$.  In particular if $w \in W(A_M, A_L)$, it induces an isomorphism between $A_M$ and $A_L$.  If $W(P,Q)$ is non trivial we say that $P$ and $Q$ are $\si$-associated. Let $c(P)=\sum_{Q\in \mm{P}} Card \> W(P, Q)$. 
  
  The proof of the  anolog of l.c. Theorem 14.3.1 (Tiling property of scattering morphisms) is entirely similar. Then one proves   the following theorem like Theorem 7.3.1 of l.c. is proved in section 14 of l.c.. Notice that one needs for this proof to establish part of this Theorem for spaces $X_P$, but this works like for $X$. 
 We recall that  $i_P$ is the map $i_{G,P}$. 
 \begin{theo}  \label{theoscat} (Scattering Theorem) Let $P=MU$, $Q=LV$, $R$ three standard $\si$-parabolic subgroups of $G$.
\\  (i) If $P$ and $Q$ are not $\si$-associated, $(i_Q)^t\circ i_P=0$.\\
(ii) If $P$ and $Q$ are  $\si$-associated, there exist $A_M\times G$-equivariant isometries 
$$S_w: L^2(X_P) \to L^2(X_Q), w\in W(P,Q)$$
where $A_M$ acts on $L^2(X_Q)$ via the isomorphism $A_M \to A_L$ induced by $w$, with the following properties:
$$ i_Q\circ  S_w= i_P, $$
$$ S_{w'}\circ S_{w}= S_{w'w}, w\in W(P, Q), w' \in W(Q,R),$$
$$ (i_Q)^t \circ i_P= \sum _{w\in W(P,Q)} S_w.$$
Let us denote by $(i_{P})^t_{disc} $  the composition of $(i_P)^t$ with the orthogonal projection to the discrete spectrum.
Finally the map
$$\sum_{P\in \mm{P}} \frac{(i_{P})^t_{disc}}{c(P)^{1/2} }: L^2(X)  \to \oplus_ {P\in \mm{P}}L^2(X_P)_{disc} $$
 is an isometric isomorphism  onto the subspaces of vectors $(f_P)_{P\in \mm{P}}\in  \oplus_{P \in \mm{P}}L^2(X_P)_{disc} $
 satisfying:
 $$S_w f_P = f_Q, w \in W(P,Q). $$ \end{theo}
 In the next section we will eplicit the maps $i_P$. 
\section{Explicit Plancherel formula}
\setcounter{equation}{0}
\subsection{Injectivity of the map  $\a' / W(A)\to \tilde{\a}' / W(\tilde{A})$}
\begin{lem}
(i) Let $A$ be a maximal $\si$-split torus and let $\tilde{A}$ be a maximal  split torus containing $A$. It  is  $\si$-stable (cf. [HH], Lemma 1.9).\\(ii) The set of non zero weights of $A$ (resp.,  $\tilde{A}$) in the Lie algebra of $G$ is a root system $\Delta(A)$ ( resp., $\Delta(\tilde{A})$) which appears as a subset of $ \a' $(resp., $\tilde{\a}'$).
\\ The set $\Delta(A)$ is equal to the set of non zero restrictions of elements  $\Delta(\tilde{A})$.\\
(iii) Let $W(A)$ (resp. $W(\tilde{A})$) be the quotient of the normalizer of $A$  (resp., $\tilde{A}$), $N_G(A)$ (resp. $N_G(\tilde{A}))$, by its centralizer, $C_G(A)$ (resp., $C_G(\tilde{A})$). 
\\ Then  $W(A)$ (resp., $W(\tilde{A})$)  identifies with the Weyl group of  $\Delta(A)$ (resp., $\Delta (\tilde{A})) $ and  is the set of restrictions  to $\a$ of the elements of $W(\tilde{A})$ which normalizes $ \a$.\\
(iv) Let $ \mu, \nu\in \a'$ which are conjugate by an element  of $W(\tilde{A})$, then they are conjugate by an element of $W(A)$.  
\end{lem}
\dem (i) follows from  [HH], Lemma 2.4.
\\ (ii) and (iii) follows from [HW], Propositions 5.3 and 5.9.
\\(iv) It is clear that one may replace $\mu$ and $\nu$ by  a conjugate  element by $W(A)$. Hence one may assume that $\mu$ and $\nu$ are dominant for some choice of  a set positive roots of $\Delta(A)$, $\Delta^+(A)$. Then we choose a  set of positive roots for $\Delta^+(\tilde{A})$ whose non zero restrictions are precisely the elements of $\Delta^+(A)$. Hence $ \mu$ and $\nu$ are dominant for $\Delta^+(\tilde{A})$ and conjugate by an element of $W(\tilde{A})$. Hence they are equal, which proves (iv).\qed 
\begin{rem} It follows from (iv) of the previous lemma that the map $\a' / W(A)\to \tilde{\a}' / W(\tilde{A})$ is injective. This allows to apply the analog of Lemma 14.2.2 of [SV].
 \end{rem}
\subsection{Coinvariants}
 Let $P=MU$ be a semistandard  $\si$-parabolic subgroup of $G$.  Let us prove:
\ber \label{sdual} Using our $G$-invariant measure on $X_P$, the smooth dual of $C^\infty _c(X_P)$ is isomorphic to   $C^\infty (X_P)$.\eer  An element of the smooth dual  of $C^\infty _c(X_P)$  is fixed by  some  compact open subgroup $J$ of $G$ and is the composition  of the $J$-average with a linear form on the  space of $J$-fixed elements of $C^\infty _c(X_P)$. A  basis of this later space  is given by  the characteristic functions  of $J$-cosets. Hence a linear form on this space  is  given by  integration of a  $J$-fixed element of  $C^\infty (X_P)$.  This proves (\ref{sdual}). 

Similarly one has:
\ber  \label{dualccxm} Using our choice of an $M$-invariant measure on $X_M$ (cf. (\ref{dxm})),  we will identify the smooth dual of $C_c^\infty (X_M)$ with $C^\infty(X_M)$. This identification depends on our choice of $\mm{X}^G_M$. \eer 
Let $(\pi,V)$ be a smooth representation of  $G$ of finite length. Let us define the space of coinvariants: \beq \label{coinv} C_c^\infty(X_P)_\pi: = Hom_\C (Hom_G(C_c^\infty(X_P), \pi), \pi).\eeq
As $Hom_G(C_c^\infty(X_P), \pi)$ is finite dimensional (cf. (\ref{finiteinv})), one has:
$$Hom_G (C_c^\infty(X_P)_\pi, \pi)= Hom_G (C_c^\infty(X_P), \pi).$$
\begin{defi} \label{canquo} If $\pi$ is a smooth admissible representation of $G$, there is a canonical projection 
   $$C_c^\infty (X_P) \to C_c^\infty (X_P)_\pi\to 0.$$ If $\pi=i^G_{P^-}\delta$, we denote this map $ i^t_{P, \delta}$
    \end{defi}
The canonical map from $C_c^\infty(X_P)$ to $C_c^\infty(X_P)_\pi$ is  defined as follows. If $f\in C_c^\infty(X_P)$, one defines $\phi  \in C_c^\infty(X_P)_\pi$ by associating to each $T \in 
Hom(C_c^\infty(X_P), \pi)$, the element $\phi(T):=T(f)$ of the space of $\pi$. It is easy to see that this map is surjective. 

Let $(\delta, E) $ be a unitary irreducible smooth representation of $M$.  Let $T \in Hom_M(C_c^\infty (X_M), \delta) $. Due to (\ref{dualccxm}),  the transpose map  $T^t $ might be viewed as an element $\tilde{T^t}$ of $ Hom (\cc{\delta}, C^\infty(X_P))$.   Let us define  $\eta_T = (\eta_{T,x})_{x\in \mm{X}^G_M}\in \mm{V}(\cc{\delta},H)$ ( cf. (\ref{vdel}) for the notation) by:
\beq \label{etat}\eta_{T,x} ( \cc{e}) := \tilde{T^t}(\cc{e}) (x), \cc{e} \in \cc{E}.Ê\eeq 
One defines $Hom_M(C_c^\infty (X_M), \delta)^{disc}$ as the space  of    $T\in Hom_M(C_c^\infty (X_M), \delta)$ such that  the image of $\tilde{T^t}$ is a discrete series for $X_M$. 
 Let us define:
 \beq \label{ccxpd} C_c ^\infty (X_P)_\delta := (Hom_M (C_c^\infty (X_M) , \delta )^{disc} )'\otimes i^G_{P^- } \delta .\eeq  
 \beq  \label{ccxdd}C_c^\infty(X_P)_\delta[\delta]=(Hom_M (C_c^\infty (X_M) , \delta )^{disc} )'\otimes \delta .\eeq 
 Hence we have:
 \beq \label{ccdd}  C_c ^\infty (X_P)_\delta= i^G_{P^-}  C_c^\infty(X_P)_\delta[\delta]. \eeq 
  It can be viewed as as  quotient of $C_c^\infty (X_P)$ as follows (cf. [ SV] before Equation   (15.12)).   From the  Lemma \ref{indcc}, one has an injective map defined by induction:
 $$0\to Hom_M ( C_c^\infty (X_M), \delta)_{disc}  \to Hom_G (C_c^\infty (X_P), i^G_{P^-} \delta).$$
 Hence, using the transpose  map and taking into account the notation (\ref{coinv}) one has a surjective map: 
 $$ C_c^{\infty} (X_P)_{i^G_{P^-}\delta } = Hom_G (C_c^\infty (X_P, i^G_{P^-} \delta)' \otimes i^G_{P^-} \delta) \to C_c ^\infty (X_P)_\delta\to 0.  $$
 Together with Definition \ref{canquo}, this shows that  \ber Ê\label{cxmquo} $C_c ^\infty (X_P)_\delta $ is a quotient of $C_c^\infty (X_P)$.\eer The smooth  dual of $C_c ^\infty (X_P)_\delta $ is denoted  $C^\infty(X_P)^{\cc{\delta}} $
 and one has $$C^\infty(X_P)^{\cc{\delta}} = Hom_M (C_c^\infty (X_M) , \delta )^{disc}  \otimes i^G_{P^-}\cc{\delta}.$$  From (\ref{cxmquo}) it can be viewed as a    subspace of  $C^\infty (X_P)$. 
\subsection{  Eisenstein integral maps  and their transpose }

\begin{defi}   We use the fact that the Eisenstein integral  associated to ${\delta_\chi}$ are well defined for  $\chi$ in the complementary set of the zero set of a non zero polynomial function on $ X(M)_\si$. For such a $\chi$, we  define a map called Eisenstein integral map: $$E_{P, \delta_\chi} \in Hom_G (Hom_M(C_c^\infty (X_M), \cc{\delta_\chi})^{disc} \otimes i^G_{P} \delta_\chi , C^\infty (X)) .$$
by $$ E_{P, \delta_\chi} (T \otimes v)= E(P, \delta_\chi, \eta_T, v), T \in  Hom_M(C_c^\infty (X_M), \cc{\delta_\chi})^{disc},  v \in i^G_P \delta_\chi.$$ 
\end{defi} 
 
We keep the notation of the preceding subsection. 
Let us denote by $ev_1$ the map 
$$ ev_1: (Hom_M(C_c^\infty (X_M), \cc{\delta_\chi})^{disc})' \otimes i^G_P \cc{\delta_\chi }  \to (Hom_M (C_c^\infty (X_M) ,\cc{ \delta_\chi } )^{disc} )' \otimes \cc{E}$$ 
defined by:
$$ev_1( \theta\otimes v)= \theta \otimes v(1),   \theta \in  (Hom_M (C_c^\infty (X_M) , \cc{\delta} )^{disc} )' , v \in i^G_{P} \cc{\delta_\chi}.$$
If $\phi \in C_c^\infty(X_M)$, let $ q_\delta( \phi)$  be the canonical  element  of $(Hom_M (C_c^\infty (X_M) , \cc{\delta} )^{disc} )' \otimes \cc{E} $ defined as follows. The later space appears as the smooth dual of $Hom_M (C_c^\infty (X_M) , \cc{\delta} )^{disc})\otimes E$  and we define $$\langle q_\delta(\phi) , T \otimes {e} \rangle := \langle {e}, T(\phi) \rangle, e \in E,  T \in Hom_M (C_c^\infty (X_M) , \cc{\delta} )^{disc}.$$  Identifying the smooth dual of  $C_c^\infty(X_M)$ to $C^\infty(X_M)$  (cf. (\ref{dualccxm})), one has also:
\beq \langle \label{qphi}  q_\delta(\phi) , T \otimes {e} \rangle = \langle \tilde{T}^t(e),\phi \rangle. \eeq  
 Let us denote, by abuse of notation,  the restriction of the transpose map of $E_{P, \delta_\chi }$ to $C_c^\infty (X) $ by 
 $E_{P, \delta_\chi }^t$ .
\begin{lem} \label{fpegal} 
 One has  $$E_{P, \delta_\chi }^t \in Hom _G (C_c^\infty (X) , (Hom_M(C_c^\infty (X_M), \cc{\delta_\chi})^{disc})' \otimes i^G_P \cc{\delta_\chi })$$
and  $$ ev_1 ((E_{P, \delta_\chi })^t ( f ))= q_\delta(f^P), f \in C_c^\infty (X) .$$
\end{lem}
\dem 
Let $e \in E$, $T \in Hom_M (C_c^\infty (X_M) , \cc{\delta} )^{disc}$. Let  $J$ be  a compact open subgroup of $G$  with a $\si$-factorization  for $P$ and such that $J_M$ fixes $e$ and $f$  and let $v_\chi:=v^{P,J}_{e, \delta_\chi}$ the element of $i^G_P \delta_\chi$ which is  invariant by $J$, whose support is equal to  $PJ$ and whose value at 1 is equal to $e$ (for the existence see e.g. [CD] Equation (3.2)). Notice that, from (\ref{sifac}),  one has:
\ber \label{suppv} $v_\chi$ has its support equal  to $PJ_{U^-}=PJ_H\subset PH$ \eer
 We will compute in two ways: 
$$I:= \langle E_{P, \delta_\chi }^t ( f), T \otimes v_\chi\rangle$$
We take into account the expression of the duality of $i^G_P \delta$ and $i^G_P \cc{\delta} $ (cf. (\ref{duali}) and (\ref{suppv})). This leads to our first expression of $I$: 
\beq \label{Iegal} I= vol(J_{U^-}) \langle ev_1( E_{P, \delta_\chi }^t ( f) ), T \otimes e  \rangle \eeq   
In order to compute $I$ in an other way  we use a transposition:
$$I= \int_{HÊ\bb G } f(\dot{g} )E_{P, \delta_\chi} (T \otimes v_\chi) (\dot{g} ) d\dot{g} .$$
 For $Re \chi$ sufficiently $P$-dominant, one has from (\ref{achtung}) and the definition of $\eta_T$ (cf.(\ref{etat})):
$$I= \int_{HÊ\bb G}  f(\dot{g} ) \sum_{x \in \mm{X}^G_M } \int_{M\cap x^{-1}.  H\bb x^{-1}. H} \tilde{T^t} (v_\chi  ( yx^{-1} g) )(x ) dy   d \dot{g}.$$
 One makes the change of variable $g'= x^{-1}.g$  and then the Fubini theorem that one can use because $f$ is  compactly supported. One gets: 
$$I=  \sum_{x\in \mm{X}^G_M } \int_{(MÊ\cap  x^{-1}. H) \bb G} f(x.g ) \tilde{T^t} (v_\chi  ( g x^{-1} ) )(x )   d \dot{g}.$$ 
We make    the change of variable $g''= gx^{-1}$.  
We use the integration formula (\ref{intumu}) and our choice of measure on $M\cap x^{-1}.H \bb M$. As $v_\chi$ has its support in $PJ_{U^-}$ and $f$ and $v_\chi$ are   $J$-invariant, one gets: 
$$I=vol (J_{U^-} ) \sum_{x \in \mm{X}^G_M }\int_{M\cap x^{-1} H \bb M }  \delta_P(m^{-1} )\int_U f(xum)du     \tilde{T^t}(v_\chi  ( m  ) )(x )  dm.$$
But the change variable $u'=m^{-1}um$ shows that:
$$I=vol (J_{U^-} ) \sum_{x \in \mm{X}^G_M }\int_{M\cap x^{-1} H \bb M }  \int_U f(xmu)du     \tilde{T^t}(v_\chi  ( m  ) )(x )  dm.$$
From the intertwining property of $T$ one has:
$$\tilde{T^t}t (v_\chi (m)) (x )) = \delta_P^{1/2}(m) \tilde{T^t}(e) ( xm).$$ With our choices of measures one deduces: 
$$I= vol (J_{U^-} ) \sum_{x \in \mm{X}^G_M } \int_{M\cap x^{-1} H \bb M } f^P (\dot{x}m ) \tilde{T^t}(e) (\dot{x}m ) dm$$ 
In other words $$I= vol(J_{U^-})\langle f^P,  \tilde{T^t}(e)\rangle, $$ 
and (\ref{qphi}) implies:
$$I= vol(J_{U^-}) )\langle q_\delta (f^P), T \otimes e\rangle.  $$
From (\ref{Iegal}) and Lemma \ref{fpegal}  (i) one deduces the equality:  $$  ev_1 ((E_{P, \delta_\chi } )^t( f )) = q_\delta(f^P).$$ \qed 
\subsection{Canonical quotient and the small Mackey restriction}
We follow the terminology of [SV], section 15. 
Let $\tau$ be  a finite length smooth  representation of $M$. If the intertwining integral:
$$A(P, P^-, \tau) : i^G_{P^-} \tau \to i^G_{P} \tau$$ 
is well defined, the {\em canonical quotient}   is the  composition:$$(i^G_{P^-} \tau)_P  \opto{{j_P(A(P, P^-, \tau))}} (i^G_P \tau)_P \to \tau$$  
where the right  map is the evaluation at $1$ (cf. [SV] Equation (15.8)).
If $\tau =  C_c^\infty(X_P)_\delta[\delta]$, 
 the canonical quotient in this case is denoted $c_\delta$  and taking into account (\ref{ccdd}) one has:
$$c_\delta: (C_c^\infty(X_P)_\delta)_P\to C_c^\infty(X_P)_\delta[\delta].$$
Let $(\pi,V)$ be a smooth representation of $G$. The Mackey restriction (cf. [SV] section 15.4.3) is the map 
$$Mack: Hom_G (C_c^\infty (X), \pi) \to Hom_M(C_c^{\infty } (X_M), \pi_P)$$
obtained by taking the Jacquet functor  to  any element $T$ of $Hom_G (C_c^\infty (X), \pi)$,  and restricting it to  $C_c^{\infty } (X_M)$ which is identied by $m_P^X$ (cf Lemma \ref{fpegal} (iii)   with a subspace of the normalized Jacquet module of $C_c^\infty (X)$. 

If $\pi= i^G_{P^-} \tau$,  and  the intertwining integral $A: i^G_{P^-} \tau \to i^G_{P} \tau$ is bijective 
the small Mackey restriction  is the composition of the canonical quotient   with the Mackey restriction $Mack$: 
$$sMack:Hom_G(C_c^\infty(X), \pi   ) \to Hom _M (C_c^\infty (X_M), \tau)$$
If  $\pi = C_c^\infty(X_P)_\delta$,  we denote by  $sMack_\delta$ the small Mackey restriction. 
If $T \in Hom_G(C_c^\infty(X), \pi) $, $$sMack_\delta (T)\in Hom_M (C_c^\infty (X_M), Hom_M (C_c^\infty (X_M) , \delta )^{disc} )'\otimes \delta)$$
\subsection{Normalized Eisenstein integrals}
\begin{defi}
Let $P$ be a semistandard $\si$-parabolic subgroup of $G$. 
We define the normalized integral
$$E^0_{P, \delta_\chi }\in Hom_G (Hom_M(C_c^\infty (X_M), \delta_\chi)^{disc} \otimes i^G_{P^-} \delta_\chi), C^\infty(X))$$ by:
$$E^0_{P, \delta_\chi}:= E_{P, \delta_\chi} \circ ( Id  \otimes A(P^-, P, \delta_\chi)^{-1})$$
which is rational in $\chi\in X(M)_\si$.
\end{defi}
By the formula of the transpose of intertwining integrals (cf. [W] IV.1(11)  and denoting by $(E^{0}_{P, \delta_\chi})^t$ the restriction of the transpose of $E^{0}_{P, \delta_\chi}$ to $C_c^\infty (X)$, one has 
$$(E^{0}_{P, \delta_\chi})^t=(Id \otimes A(P, P^-, \delta_\chi)^{-1}) \circ( E_{P, \delta_\chi})^t. $$
From this it follows 
$$ sMack ((E^{0}_{P, \delta_\chi})^t) \in Hom (C_c^\infty (X_M), Hom_M (C_c^\infty (X_M) , \delta_\chi )^{disc} )'\otimes \delta_\chi)$$ 
 is equal to   $$ev_1 j_P (A(P, P^-, \delta)\circ(E^{0}_{P, \delta})^t).$$ 
From  Lemma  \ref{fpegal}, one deduces: \ber \label{smei}The map  $sMack ((E^{0}_{P, \delta_\chi})^t)$  is equal to the map $q_{\delta_\chi}$.
 \eer 
 \ber  \label{smr} Our definition of normalized Eisenstein integrals differs from the one in [SV], Equation (15.30) for $G$-split and $X$ spherical. Here we do not use the Radon transform, but we use that the opposite of a $\si$-parabolic subgroup is a $\si$-parabolic subgroup.  From (\ref{smei}), our  Eisenstein integrals maps  have   the same small Mackey restrictions restrictions than the ones  defined  in l.c.    (cf. (15.36)). \eer 
  \subsection{Explicit Plancherel formula}
  Let $$L^2(X_M)_{disc} =  \int^{\oplus} _{ \hat{M}} \cc{I} _\delta d\nu(\delta) $$
where $\cc{I}$ is a unitary representation of $M$ isomorphic to a direct sum of copies of $\delta$.  
From Lemma  \ref{inddisc} , one has  $$L^2(X_P)_{disc}= \int^{\oplus} _{ \hat{M}}  \cc{H} _\delta d\nu(\delta) $$  
  where $\cc{H}_\delta $ is the unitarily induced representation from $P^-$ to $G$ of $ \cc{I}_\delta$.   Let $ \cc{H}_\delta^\infty $ be its space of smooth vectors. With the notation of (\ref{ccxpd}), its space of smooth vectors is equal to $C^\infty (X_P)^{\delta}$.

  Let  $f \in C_c^\infty (X_P)$  and let us write  its discrete component $$f_{disc} = \int_{\hat{M} } f^\delta d\nu_{disc} (\delta),  $$ where $f^\delta \in C^\infty (X_P)^\delta$.
 Its image by the Bernstein morphism $i_P (f)$ satisfies:
 $$i_P (f_{disc} )= \int_{\hat{M}} i_{P, \delta}  (f) d \nu( \delta).   $$
 for some   maps  $i_{P, \delta}:\cc{H}_\delta^\infty    \to C^\infty (X)$ defined for almost all $\delta$ (cf. [SV] Equation (15.6)).

  One has the analog of Lemma 15.4.4 of [SV]. 
   As the analogous of section 15.6 of [SV] is identical, together with (\ref{smei}),  this leads to the  analog  of Th 15.5.5 in [SV]:
  \begin{prop} The small Mackey restriction of $ sMack ((E^{0}_{P, \delta})^t)$ and $i^t_{P, \delta_\chi}$ are equal for almost all $\chi\in X(M)_{\si, u} $. 
  \end{prop} 
  
  Also by   the uniqueness  result of [BD] recalled in (\ref{uniqueext}), for almost all $\chi\in X(M)_{\si, u}$, every element $F$ of $Hom_G (C_c^\infty(X), i^G_{P^-}\delta_\chi)$ is given in term of the normalized Eisenstein integral i.e. is of the form
  $$F=E(P, \delta_\chi, \eta_T, v)\circ  A(P^-,P, \delta_\chi)^{-1}$$ for a unique $ T \in  Hom_M(C_c^\infty (X_M), \cc{\delta_\chi})$. Using (\ref{smei}) or rather its immediate generalization by replacing  $Hom_M(C_c^\infty (X_M), \cc{\delta_\chi})^{disc} $ by $ Hom_M(C_c^\infty (X_M), \cc{\delta_\chi})$ one sees that    the small Mackey restriction of $F$ is equal to $T$. Hence one has: \begin{prop}
  The small Mackey restriction  $$sMack:Hom_G(C_c^\infty(X), i^G_{P^-} \delta_\chi   ) \to Hom _M (C_c^\infty (X_M), \delta_\chi )$$is injective for almost all $\chi\in X(M)_{\si, u} $. 
   \end{prop}
   
   \begin{cor} For almost all $\chi\ in X(M)_{\si, u} $, one has:
   $$ i_{P, \delta_\chi}= E^{0}_{P, \delta_\chi}.$$
   \end{cor}

 \begin{theo}
 Let  $f \in C_c^\infty (X_P)$ and  let us write  its discrete component $$f_{disc} = \int_{\hat{M} }^\oplus f^\delta d\nu_{disc} (\delta),  $$ where $f^\delta \in C^\infty (X_P)^\delta$.\\
 Its image by the Bernstein morphism $i_P (f)$ satisfies:
 $$i_P (f)(x)= \int_{\hat{M}} E^0_{P, \delta}( f^\delta)(x)d \nu( \delta), x \in X.$$
 \end{theo}
 In combination with the scattering theorem (cf. Theorem \ref{theoscat}), one deduces: 
 \begin{theo}  The norm on  $L^2(X)_P$, $\Vert . \Vert_P$,   admits  the decomposition:
 $$Ê\Vert \Phi\Vert _P^2 = \frac{4}{ Card ( W(A_P, A_P))}\int_{\hat{M} }  \Vert  E^{ 0t} _{P, \delta} (\Phi)  \Vert ^2_\delta d\nu(\delta),$$
where the measure and norms on the right hand side of the equality are the discrete part of the Plancherel decomposition of $L^2(X_P) $.
  \end{theo}

\section{Appendix: Rational representations} \label{appendix}
\setcounter{equation}{0}
In this section we establish some results on rational representations of $G$ which are needed to extend the results of [L] and [BD], which are established when $\F$ is of characteristic zero,   to the case where $\F$ is simply of characteristic different from 2.
\subsection{Rational representations and parabolic subgroups}
Let $\underline{G}$ be  a reductive  algebraic group defined over  a non archimedean local field $\F$, whose group of $\F$-points is equal to $G$. We will use similar notations for the subgroups of $G$. \\
Let $A_0$ be a maximal split torus of $G$ and  let $P_0$ be  a minimal parabolic subgroup of $G$. Let $P$ be a parabolic subgroup of $G$ which contains $P_0$. Let $T$ be a maximal $\F$-torus  of $\underline{G}$ which contains $\underline{A}_0$. Let $B$ be  a Borel subgroup of  $\underline{G}$, which contains  $T$ and contained in the opposite parabolic subgroup to  $P_0$ which contains $A_0$,  $\underline{P}^-_0$.  One denotes by  $\Sigma(T)$ the set of roots of $T$ in the Lie algebra of  $\underline{G}$. One denotes by  $\Lambda(T)$ (resp.,   $\Lambda(T)_{rac}$) the weight lattice  (resp., the root lattice) of  $T$  with respect to   $\underline{G}$. We adopt similar notations for $A_0$. Let  $\Gamma$ be the absolute   Galois group  of  $\F$ which acts on these lattices. Let  $\Lambda ^+(T) $  be the set   of dominant weights for $T$ relative to  $B$.  Let  $\LL^+(A_0)$ (resp., $\LL^+(A_0)_{rac}$) the set of dominant elements for $P^-_0$ of $\LL^+(A_0) $ (resp., $\LL(A_0)_{rac}$).  \begin{defi} One  denotes by  $\Lambda^+_M(T)$ the set of elements   $\lambda$  of  $\Lambda^+(T)$such that  $G$ has  a rational finite dimensional  irreducible representation, defined over $\F$,   with highest weight $\lambda$ relative to $B$,  $(\pi_\l, V_\l)$,  with the following property:
\ber Any non zero vector of weight $\l$  under $T$, $v_\l$,  
 transforms under $M$ by a rational character of $M$, denoted $\LL$. \eer 
The subset of $\l\in \LL^+_M(T)$   which satisfies the following property is denoted by $\LL^{++}_M(T)$.
\ber There exists  exists $v'_\l$ in the dual $  V'_\l$  of $V_\l$, invariant by  $U$ and such that the  coefficient $c_\l (g)= \vert \langle \pi_\l(g)v_\l, v'_\l\rangle \vert_\F$ is equal to zero on the complementary subset in  $G$ of  $UMU^-$.  \eer 
\end{defi}
The goal of this subsection is to produce sufficiently many elements of $\LL^+_M(T)$.\\
\begin{prop} \label{nl}
(i) Let $T_{an} $ be the anisotropic component of $T$. There exists $n\in \N^*$ such that any element $\l$ of $n\LL^+(A_0)$ extends uniquely to  an element $\mu
 $ of $\LL(T)_{rac} $ trivial  on $T_{an}$.\\
(ii) If $\l$ is  orthogonal to the simple roots  of  $A_0$ in  the Lie algebra of   $U_0^-\cap M$ then $\mu$ is element of $\LL^+_M(T)$.\\
(iii) If moreover $\l$ is not  orthogonal to the other  simple roots of $A_0$ in the Lie algebra of $U_0^-$,  $\mu $ is an element of $\LL^{++}_M(T)$. 
\end{prop}
For the proof we need several lemmas.

Let $\beta$ be an element of the set, $\Sigma(A_0)$, of roots  of $A_0$ in the Lie algebra of $G$.  One defines:  $$\underline{\beta}: = \sum_{ \aa \in \Sigma(T), \aa_{\vert A_0}= \beta} \aa. $$
One sees easily that: \ber There exists $n'\in \N^*$ such that, for all $\beta \in \Sigma(A_0)$,  there exists  $n'_\beta\in \N^* $such that  $ n'_\beta\underline{\beta}_{\vert A_0}= n' \beta$. \eer
We fix, once for all, such  integers $n'$ and $n'_\beta$
\begin{lem} \label{extend}
 Every element $\l$ of $n'\LL_{rac}(A_0)$ extends uniquely to an element $\mu $ of $\LL_{rac}(T)$ trivial on the anisotropic component $T_{an}$ of $T$,  invariant by $\Gamma$ and by $W(\underline{M}_0, T)$. 
\end{lem} 
Let us denote by  $(n'\Lambda_{rac}(A_0))^{\tilde{}}$ the lattice generated by the  $n'_\beta \underline {\beta}$, $\beta \in \Sigma(A_0)$. From their definition, one sees that the elements of  $(n'\Lambda_{rac}(A_0))^{\tilde{}}$ are invariant under  $\Gamma$ and are elements of  $\Lambda_{rac}(T)$. One remarks that every element $\mu$  of   $(n'\Lambda_{rac}(A_0))^{\tilde{}}$ is  invariant by the  Weyl  group  of $\underline{M}_0$ relative to $T$, $W(\underline{M}_0 , T)$.\\
    Let us show any  element  $\mu$ is trivial on $T_{an}$.  One can choose $T$ such that it contains a maximal torus defined over $\F$, $T_1$, of the derived group of $\underline{M}_0$.  Actually, by conjugacy, one sees that any $T$ has this property. Moreover $T$ contains  the maximal anisotropic  torus  $C_{an}$ of the center of $M_0$.    The product  
  $ T_1 C_{an} \underline{A}_0$ is a torus. For reasons of dimension it is a maximal torus of $G$. Hence $T= T_1 C_{an} \underline{A}_0$. Notice that $T_1 C_{an}$ is the anisotropic component $T_{an}$  of $T$. As  $\mu $ is  $W(\underline{M}_0 , T)$-invariant, the restriction of  $\mu $ to  $T_1$ is trivial. As  $C_{an}$  is anisotropic, the invariance by $\Gamma$of  $\mu$ shows that its restriction to   ˆ $C_{an}$ is trivial.  This proves the existence part of the Lemma. As $T= T_{an} \underline{A}_0$ the unicity follows. \qed 

\begin{lem} \label{ratmu}
(i) There exists $n\in \N$ such that $n\LL(A_0) \subset n' \LL_{rac}(A_0)$.\\
(ii) If $ \l\in n\LL^+(A_0)$, its extension $\mu$ to $T$ given by the preceding lemma is the highest weight of a rational representation of $G$,   defined over $\F$, denoted $(\pi_{\mu}, V_{\mu} ).$
\end{lem}
\dem 
 (i) The lattice 
 $n'\LL _{rac}(A_0)$ is contained in the lattice $   \LL(A_0)$. As these lattices are of  the same rank,  there exists  $n\in \N^*$ such that 
  $ n\LL(A_0) \subset n'\Lambda_{rac}(A_0).$ 
  \\ (ii) From (i) and the preceding lemma, if  $ \l\in n\LL^+(A_0)$, $\mu$ is in  $\LL_{rac}(T) \subset \LL(T)$. Moreover
  if $\aa$ is a root  of  $T$ in the Lie algebra of $\underline{G}$, $\langle\mu,\aa \rangle =\langle \l,\aa_{\vert A_0} \rangle$. Hence $\mu$ is a dominant weight. From the preceding Lemma, it is invariant by $\Gamma$. Then  [T],  Theorem 3.3 and Lemma 3.2 implies (ii). \qed 
  \begin{lem}
Let  $\l \in n \LL^+(A_0)$ and $\mu$ as  in Lemma \ref{extend}. Then, with the notation of the preceding lemma,  $M_0$ acts on  a non zero highest weight vector of $(\pi_\mu, V_\mu)$ by a rational character of $M_0$ again denoted  by $\mu$.
\end{lem}
\dem 
 As $\pi_\mu$ is defined over $\F$, it is enough to prove that $v_\mu$ transforms under a rational character of $M_0$.  In order to prove this, one can work with the algebraic closure. The invariance of $\mu$ by $W(\underline{M}_0 , T)$ (cf. Lemma \ref{extend}), the fact that the space of weight  $\mu$ in $V_\mu$ is of dimension one (cf. [Hu], Proposition 33.2) together with the  Bruhat decomposition of  $\underline{M}_0$ allow to prove the Lemma. \qed 
 \\ {\em  Proof of Proposition \ref{nl}}\\
 Let $\l \in n \LL^+{(A_0)}$  be as in the statement of Proposition \ref{nl} (ii) i.e. $\l$ is  orthogonal to the simple roots  of  $A_0$ in  the Lie algebra of   $U_0^-\cap M$.
 Let $\mu$ be as in Lemma \ref{extend}. Let $(\pi_\mu, V_\mu) $ be as in Lemma \ref{ratmu}, and let $v_\mu$ be a non zero highest weight vector.  One has to prove   that  $v_\mu$ transforms under $M$ by a rational character of  $M$ that we will still denote by $\mu$. It is enough to prove that the line $\F\mu$ is stable by the action of $M$.   One shows,  using the preceding Lemma, by a proof analogous to the one of  [Hu],   Proposition  31.2 and using the density  of  $U_0^-M_0U_0$ in  $G$,  that the $A_0$-weight space  of $V_\mu$ for the weight $\l$ is one dimensional. The Weyl group, $W(M, A_0)$, of $M$ relative to $A_0$ fixes $\l$  from the hypothesis on $\l$. One finishes the proof of our assertion on the action  $M$ on $v_\mu$ by using the Bruhat decomposition of  $M $ relative to $P_0^- \cap M$.  Hence $\mu\in \LL^+_M(T)$.
 
 Now we assume moreover that  $\l$ is not  orthogonal to the other  simple roots of $A_0$ in the Lie algebra of $U_0^-$,  $\mu $ is an element of $\LL^{++}_M(T)$. 
  One sees like in l.c. that the weights of $A_0$ in $V_\mu$ are  of the form  $\nu = \l + \sum_{\beta \in \Delta(P_0) }c_\beta \beta$, where for all $\beta$ in the set of simple roots $\Delta(P_0)$ of $A_0$ in the Lie algebra of $U_0$,  $c_\beta\in \N$  (we recall that   $B$ is contained in   $\underline{P}^-_0$). We consider the hyperplane of  $V_\mu$  generated by the $A_0$-weight spaces for the weights distinct from $\l$.  From what we have just established on the weights of $A_0$ in $V_\mu$ and from  [Hu], Proposition 27.2, one sees that this  hyperplane    is stable by   $U$. Hence the linear form on $V_\mu$, $v'_\mu$,  which vanishes on this hyperplane and which takes the value 1 on $v_\mu $  transforms  under a rational character of the   unipotent group $U$. This implies that $v'_\mu$ is invariant by  $U$. \\
Let  $c_\mu $ be   the real valued function on $G$ defined by:
$$c_\mu(g)= \vert \langle  \pi_\mu(g)v_\mu, v'_\mu\rangle  \vert _\F, \> g \in G. $$
Let us show that our hypothesis on  $\lambda$ implies that   $c_\mu $ vanishes on the complimentary subset of $UMU^-$ in $G$.
From the Bruhat decomposition for $P_0^-$ an element of this complimentary subset  can be written  $g=u wmu^-$,  where $w$ represents an element of $ W(G, A_0)$  which is not in $W(M,A_0)$. Then $c_\mu(g)$ is proportional to  $\vert \langle 
\pi_{\mu}(w)v_\mu, v'_\mu\rangle   \vert_\F$ .  But 
$\pi_\mu(w)v_\mu$ is of weight  $w  \l$ under  $A_0$. This weight is  
distinct of  $\lambda$  as   $\l$  is not orthogonal to the simple roots of $A_0$ in the Lie algebra $U_0^-$ which are not roots of $A_0$ in $M$ . This implies that  $\vert \langle 
\pi_{\lambda}(w)v_\lambda, v'_\lambda\rangle    \vert_F$ is equal to zero.  Hence   $c_\mu(g)$ is equal to zero as wanted.  This achieves the proof of the Proposition.  \qed  

\subsection{$H$-distinguished rational representations of $G$}  \label{lagier}
Proposition \ref{nl} allows to extend the results of [BD] section 2.7 and especially Propositions 2.9, 2.11 to a non archimedean local field, $\F$, of characteristic different from 2. 
Let  $\Sigma(G, A_{0})$ (resp., $\Sigma(P_{0}, A_{0})$ or  $\Sigma(P_{0})$ the set of roots  of $A_0$ in the Lie algebra of  $G$ (resp., $P_0$).  We denote  by 
$\Delta(P_{0})$ the set of simple roots of  $\Sigma(P_{0})$. 

Let $P=MU$ be  a standard  $\si$-parabolic subgroups of $G$. We will use the notation of the main body of the article. Let us assume that  $A_\vid \subset A_0$, which is automatically $\si$-stable, and $P_0\subset P_\vid$.   Let  $\{ \aa_1, \dots , \aa_{m_0} \} $ be the simple roots of $\Sigma(P_{0})$  written in such a way that 
$\{ \aa_1, \dots , \aa_{m_\vid} \} $ are the simple roots in the Lie algebra of    $U_\vid $, 
$\{ \aa_1, \dots , \aa_{m} \} $ are the simple roots in the Lie algebra of  $U$. 
 One has   the fondamental weights of 
$\Sigma(P_{0}, A_{0})$, $\delta_1, \dots, \delta_l$. \\ Let   $i=1,\dots, m $ and  $\l_i= n\delta_i$ with $n $ as in Proposition \ref{nl}.  From this proposition,  there exists a unique  rational character of $T$, $\mu$,  trivial on $T_{an}$ and whose restriction to $A_0$ is  equal to $\l_i$ and such that $\l_i \in \LL^+_{M}(T) $ and   $\mu$  is the highest weight of an irreducible finite dimensional rational  representation of $G$, denoted by $(\pi_\mu, V_\mu)$. Moreover if $v_\mu $ is  a non zero highest weight vector in $V_\mu$,  
 the space  $\F v_\mu$ is 
$P$-invariant. We denote again by $\mu$ the rational character of $M$ which describes the action of $M$ on $v_\mu$. One  denotes by  $v'_\mu$
the unique element of 
$ V_{\mu}' $ of weight 
$\mu^{-1}$ under 
$M $ and such that   $\langle v'_\mu,v_\mu \rangle = 1$. 

Let 
$\nu:=\mu(\mu^{-1}\circ\sigma)\in \Lambda(T)$  and let  $ (\tilde{\pi}_\nu,  \tilde{V}_\nu )$ be the rational representation of $G$   $(\pi_\mu \otimes (\pi_\mu'\circ
\sigma),V_{\mu }\otimes V_{\mu}')$. Let  $\tilde{v}_\nu:=v_\mu\otimes v_\mu'$  which is of weight 
 $\nu $ under  the representation $ \tilde{\pi}_\nu$ restricted to  $M $. Then there exists a non zero  $H$-invariant vector, under $\tilde{\pi}_\nu$  in $\tilde{V}_\nu'=(V_{\mu}\otimes
V_{\mu} ')'\simeq V_{\mu}'\otimes V_{\mu} \simeq End V_\mu$, namely the identity that we will denote $e_{\nu,  H}'$. It satisfies $\langle
e_{ \nu,H}', \tilde{v}_\nu \rangle = 1$.

Let us show that $\nu= 2 \mu$.   As $\si$ preserves $T_{an}$,  the character $\mu^{-1}\circ\sigma$   of $T$ is trivial on $T_{an}$. Its  restriction to $A_0$ is equal to $\l$. From the unicity statement of $\mu$, it is equal to $\mu$. This proves our claim.

From this it follows that 
\ber \label{bdl} Proposition 2.9 and 2.11  of  [BD] extend to a non archimedean local field, $\F$,  of characteristic different from 2.
 This shows that the results of [BD],  section 2.8, 2.9 are valid for such a field .  Also, the Lemma 1 (resp., section 3.2) in [L] is true also for such a field $\F$ due to Proposition 2.3 of [CD] (resp., the Proposition \ref{nl} of the present article). Hence the results of  [L] are valid for such a field $\F$.\eer 



\section{References}

\noindent [AAG]  Aizenbud A., Avni N., Gourevitch D.,  Spherical pairs over close local fields. Comment. Math. Helv. 87 (2012) 929--962. 

\noindent {BenO]  Benoist  Y., Oh H.,  Polar decomposition for p-adic symmetric spaces. Int. Math. Res. Not. IMRN 2007, Art. ID rnm121. 

\noindent{[Ber]  Bernstein J.,  On the support of Plancherel measure.
J. Geom. Phys. 5 (1988)  663--710 (1989).

\noindent [BD] Blanc P., Delorme P.,  Vecteurs distributions H-invariants de
reprŽsentations induites pour un espace symŽtrique rŽductif $p$-adique
$G/H$, Ann. Inst. Fourier, 58 (2008),  213--261. 

\noindent [CD) Carmona J.,  Delorme P., Constant term of H-forms  arXiv:1105.5059

\noindent[D] Delorme P., Constant term of smooth $H_\psi$-invariant functions,  Trans. Amer. Math. Soc.  362  (2010),  933--955.  

  \noindent [HH] Helminck A.G., Helminck G.F., A class of parabolic $k$-subgroups
associated with symmetric
$k$-varieties.  Trans. Amer. Math. Soc.  350  (1998)
4669--4691.

\noindent[HW] Helminck A. G.,  Wang  S. P., 
On rationality properties of involutions of reductive groups. Adv. Math. 99 (1993) 26--96.

\noindent [Hu] Humphreys J. E, Linear algebraic groups, Graduate Text In Math. 21, Springer, 1981.
 
\noindent [KT1] Kato S., Takano K., Subrepresentation theorem for $p$-adic symmetric spaces,  Int. Math. Res. Not. IMRN 2008, no. 11. 

\noindent [KT2] Kato S., Takano K., Square integrability of representations on $p$-adic symmetric spaces. J. Funct. Anal. 258 (2010)1427Ð-1451.

\noindent [L]  Lagier N., Terme constant de fonctions sur un espace symŽtrique rŽductif p-adique,J. of Funct. An., 254 (2008) 1088--1145.

\noindent [R]  Renard D., ReprŽsentations des groupes rŽductifs $p$-adiques. Cours SpŽcialisŽs, 17. SociŽtŽ MathŽmatique de France, Paris, 2010.

 \noindent [T] Tits J., ReprŽsentations linŽaires d'un groupe rŽductif sur un corps quelconque,  J. Reine Angew. Math. 247 1971 196--220. 

 \noindent [SV]  Sakellaridis Y., Venkatesh A.,  Periods and harmonic analysis on spherical varieties, arXiv:1203.0039

  \noindent [ W] Waldspurger J.-L.,  La formule de Plancherel pour les groupes
$p$-adiques (d'aprs Harish-Chandra), J. Inst. Math. Jussieu  2  (2003),  235--333.
\\\\
 Aix Marseille Universit\'e\\
CNRS-IML, FRE 3529\\
163 Avenue de Luminy\\ 13288 Marseille Cedex 09\\  France.\\          
 delorme@iml.univ-mrs.fr
 \\\\
 The author has been supported by the program ANR-BLAN08-2-326851 and by the Institut Universitaire de France  during the elaboration of this work.}
 \end{document}